\documentclass[11pt]{amsart}

\usepackage{epsf,amssymb}

\newtheorem{thm}{Theorem}[section]
\newtheorem{prop}[thm]{Proposition}

\newtheorem{lemma}[thm]{Lemma}
\newtheorem{cor}[thm]{Corollary}

\newtheorem{q}[thm]{Question}

\newcommand{\R}{\mbox{\bf{R}}}
\newcommand{\Z}{\mbox{\bf{Z}}}

\newcommand{\bdry}{\partial}
\newcommand{\cal}{\mathcal}
\newcommand{\sa}{\rightsquigarrow}
\renewcommand{\sc}{\rightarrowtail}
\newcommand{\s}{\vskip.12in}
\newcommand{\n}{\noindent}

\newcommand{\be}{\begin{enumerate}}
\newcommand{\ee}{\end{enumerate}}

\begin{document}
\title{Gluing tight contact structures}

\author{Ko Honda}
\address{University of Georgia, Athens, GA 30602}
\email{honda@math.uga.edu}
\urladdr{http://www.math.uga.edu/\char126 honda}

\date{First version: November 28, 1999.  This version: July 3, 2000.}


\begin{abstract}
	We prove gluing theorems for tight contact structures.  In particular, we rederive (as special cases)
	gluing theorems due to Colin and Makar-Limanov, and present an algorithm for determining whether
	a given contact structure on a handlebody is tight.  As applications, we construct a tight contact structure
	on a genus 4 handlebody which becomes overtwisted after Legendrian $-1$ surgery, and study certain
	Legendrian surgeries on $T^3$.
\end{abstract}
\maketitle

\section{Introduction}

The goal of this paper is to provide an updated account of the theory of tight contact structures on
3-manifolds.  The emphasis is on purely 3-dimensional techniques, which arises from
a Haken decomposition theory, where the cutting manifolds are {\it convex surfaces}.
In the introductory sections of this paper, we will briefly explain the notion of a convex surface,
which appear to be the best kind of cutting surface for a decomposition, as well as the notion
of a {\it bypass}, which plays a role in flaking off $\Sigma\times I$ layers, where $\Sigma$ is a
convex surface and $I=[0,1]$.  We will explain how to cut along convex surfaces with Legendrian boundary --
in order to ensure that the boundary that we cut along is Legendrian, we use a preparation theorem
called {\it Legendrian realization}.  We are then able to perform a {\it convex decomposition} on
$M$.

The following question then arises naturally.
Suppose we want to construct a tight contact structure using the
reverse procedure.  If $M\sa M'$ is one step of the convex decomposition, and $M'$ is tight,
what are the conditions for $M$ to be tight?  We provide one answer in the form of Theorems \ref{gluing} and
\ref{classification}.  Theorem \ref{gluing} provides an explicit algorithm for determining whether
a prescribed contact structure on a handlebody is tight, and allows us, at least in theory, to classify tight
contact structures on a handlebody with prescribed boundary conditions in finite time.
Theorem \ref{classification} is a general gluing/classification theorem, which is usually rather difficult to
verify in practice.  However, as special (combinatorially simple) cases of Theorem \ref{classification},
we are able to recover gluing theorems due to Colin \cite{Co97}, \cite{Co99} and Makar-Limanov \cite{ML}
on tightness preserved under connect sum operations, and restricted boundary connect sum
operations.

We will discuss some specific examples of tight contact structures which
can be constructed using the Gluing Theorem.
(1) We will exhibit a tight contact structure on a genus 4 handlebody which
becomes overtwisted after Legendrian surgery.  This is the first example of a tight contact structure
(albeit not closed) which does not remain tight after Legendrian surgery.
(2) We will study some Legendrian surgeries of $\xi_n$, $n\in \Z^+$,
the distinct tight contact structures on
$T^3=S^1\times T^2$.

\subsection{Convex surfaces}

We assume that the reader is familiar with the basic ideas of contact topology in
dimension 3 (see for example \cite{E92}), especially the dichotomy between
{\it tight} contact structures and {\it overtwisted} contact structures.  Since
overtwisted contact structures are relatively well-understood by the work of Eliashberg \cite{E89}, our focus is on
tight contact structures.

Let $M$ be an oriented, compact 3-manifold (possibly with boundary), and let
$\xi$ be a {\it positive} contact structure which is co-oriented by a global 1-form $\alpha$ with
$\alpha\wedge d\alpha>0$.
We define {\it Legendrian curves} to be  {\it closed} curves which are everywhere tangent to $\xi$, as opposed
to {\it Legendrian arcs}.

If $X$ is a manifold and $Y$ is a submanifold, then we will use the notation $X\backslash Y$
to mean the metric closure of the complement of $Y$
in $X$.

An oriented properly embedded surface $\Sigma$ in $(M,\xi)$ is called {\it convex} if
there is a vector field $v$ transverse to $\Sigma$ whose flow preserves $\xi.$
The {\it dividing set} $\Gamma_\Sigma$ of $\Sigma$  {\it with respect to $v$} is the set of points $x$ satisfying
$v(x)\in \xi(x)$.  The isotopy type of $\Gamma_\Sigma$ is independent of the choice
of $v$ -- hence we will usually call $\Gamma_\Sigma$ {\it the dividing set} of $\Sigma$.
$\Gamma_\Sigma$ is a union of pairwise disjoint smooth curves which are transverse
to the {\it characteristic foliation} $\xi|_\Sigma$.    Denote the
number of connected components of $\Gamma_\Sigma$ by $\#\Gamma_\Sigma$.
$\Sigma\backslash\Gamma_\Sigma = R_+-R_-$, where $R_+$ is the subsurface where
the orientations of
$v$ (coming from the normal orientation of $\Sigma$) and the normal orientation of $\xi$
coincide, and $R_-$ is the subsurface where they are opposite.

We will not assume that our convex surfaces are closed.  In fact, in most applications,
our convex $\Sigma$ will be a compact convex surface with {\it Legendrian boundary}.
This will impose a condition on each connected component $\gamma$ of $\bdry\Sigma$
-- namely, the {\it twisting number} $t(\gamma,Fr_\Sigma)$ of $\gamma$ relative to the
framing $Fr_\Sigma$ induced from $\Sigma$ must be nonnegative.    Here we are using the
convention that left twists are negative.

\vskip.15in
\noindent
{\bf Key Principle:} It is the dividing set $\Gamma_\Sigma$ (not the exact characteristic
foliation) which encodes the essential
contact topology information in a neighborhood of $\Sigma$.
\vskip.15in

To make this idea more precise, we will now present Giroux's Flexibility Theorem.
If ${\cal F}$ is a singular foliation on $\Sigma$, then a disjoint union of properly embedded
curves $\Gamma$ is said to {\it divide} ${\cal F}$  if there exists some $I$-invariant
contact structure $\xi$ on
$\Sigma\times I$ such that ${\cal F}=\xi|_{\Sigma\times \{0\}}$ and $\Gamma$ is the dividing
set for $\Sigma\times \{0\}$.

\begin{thm}[Giroux \cite{Gi91}]\label{flexibility}
	Let $\Sigma$ be a convex surface with characteristic foliation $\xi|_\Sigma$,
	contact vector field $v$, and dividing set $\Gamma$. If  $\mathcal{F}$ is another
	singular foliation on $\Sigma$ divided by $\Gamma$, then there is an
	isotopy $\phi_t$, $t\in[0,1]$, called an {\em admissible isotopy} of $\Sigma$,
	such that $\phi_0(\Sigma)=\Sigma,$
	$\xi|_{\phi_1(\Sigma)}=\mathcal{F}$, the isotopy is fixed on $\Gamma$,
	and $\phi_t(\Sigma)$ is transverse to $v$ for all $t$.
\end{thm}

The following is Giroux's criterion for determining which convex surfaces have neighborhoods which
are tight:

\begin{thm}[Giroux's criterion] If $\Sigma\not = S^2$ is a convex surface in a contact manifold
$(M,\xi)$, then $\Sigma$ has a tight neighborhood if and only if $\Gamma_\Sigma$ has no homotopically
trivial curves.  If $\Sigma=S^2$, $\Sigma$ has a tight neighborhood if and only if
$\#\Gamma_\Sigma=1$.
\end{thm}

\begin{figure}
	{\epsfysize=1.5in\centerline{\epsfbox{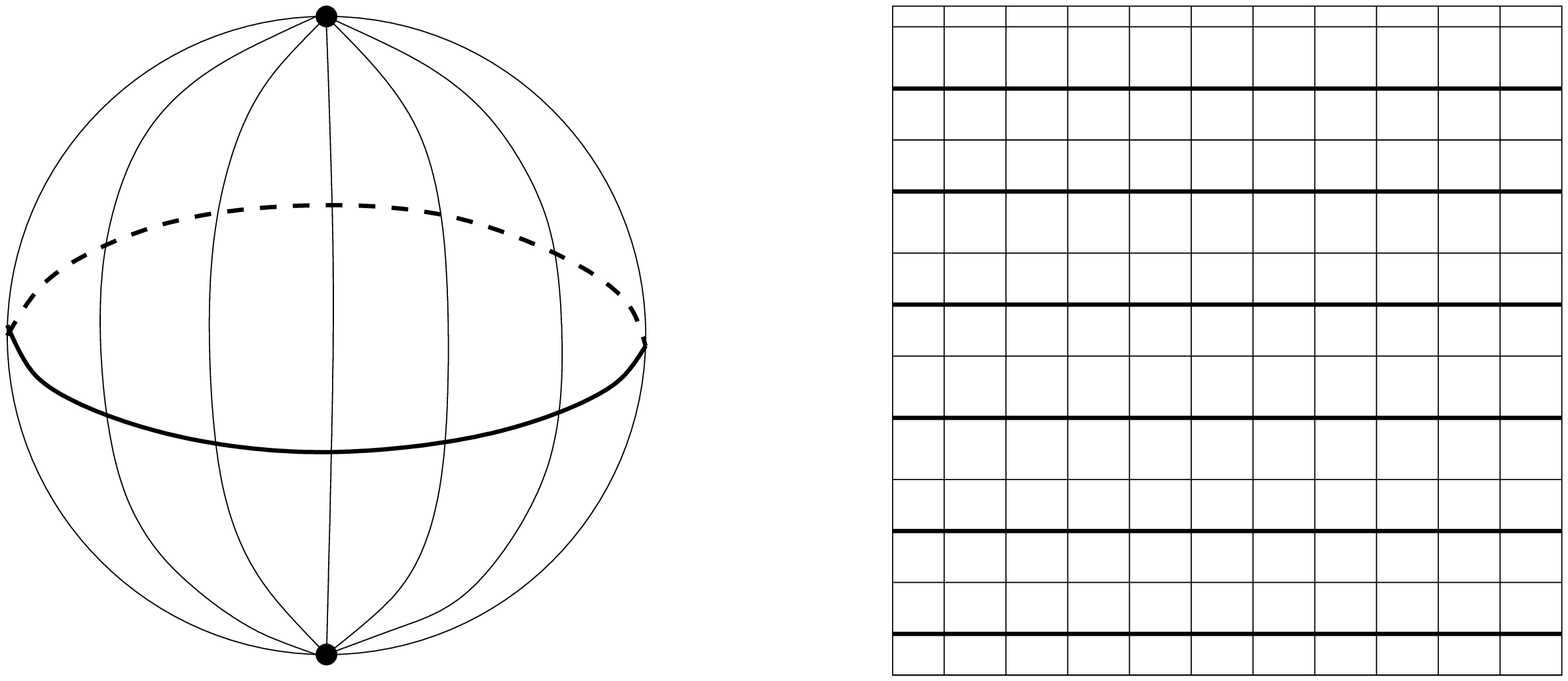}}}
	\caption{Dividing curves for $S^2$ and $T^2$}
	\label{fig1}
\end{figure}

\noindent
{\bf Examples:}  The following are some examples of convex surfaces that can exist inside
tight contact manifolds.
\be
\item $\Sigma=S^2$.  Since $\#\Gamma_\Sigma=1$, there is only one possibility.  See Figure
\ref{fig1}. Note that any time there is more than one dividing curve the contact structure is
overtwisted.   In Figure \ref{fig1}, the thicker lines are the dividing curves.

\item $\Sigma=T^2$.  Since there cannot be any homotopically trivial curves, $\Gamma_\Sigma$
consists of an even number ($>0$) of parallel homotopically essential curves. If we
identify $T^2\simeq \R^2/\Z^2$ so that the dividing curves are horizontal, we have the situation in Figure \ref{fig1}.
Note that the sides are identified and the top and bottom are identified.  The thin horizontal lines
are tangencies, called {\it Legendrian divides}, and the vertical lines form a family of Legendrian curves, called
{\it Legendrian rulings}.
\ee

\subsection{Legendrian realization principle}

In this section we present the {\it Legendrian realization principle} -- a criterion for determining
whether a given curve or a collection of curves and arcs can be made Legendrian after a
perturbation of a convex surface $\Sigma$.  The result is surprisingly strong -- we can realize
almost any curve as a Legendrian one.   Our formulation of Legendrian realization is a
generalization of Kanda's \cite{K98}.
Call a union of closed curves and arcs $C$ on a convex surface $\Sigma$ with Legendrian
boundary {\it nonisolating}
if  the following hold:
\be
\item $C$ is transverse to $\Gamma_\Sigma$.
\item Every arc in $C$ either begins and ends on
$\Gamma_\Sigma$, or begins and ends on $sing(\bdry \Sigma)$.  Here assume that $\Sigma$ has
been normalized near $\bdry \Sigma$ so that the singular set $sing(\bdry \Sigma)$ of $\bdry \Sigma$ consists
solely of half-elliptic points.
\item $C$ is pairwise disjoint, with the possible exception of
arcs having common endpoints on $sing(\bdry \Sigma)$.
\item Every component of $\Sigma\backslash (\Gamma_\Sigma \cup
C)$ has a boundary component which intersects $\Gamma_\Sigma$.
\ee

\begin{thm}[Legendrian realization] \label{lerp}
Consider $C$, a nonisolating collection of closed curves and arcs, on a convex surface $\Sigma$
with Legendrian boundary.
Then there exists an
admissible isotopy $\phi_t$, $t\in[0,1]$ so that
\be
\item $\phi_0=id$,
\item $\phi_t(\Sigma)$ are all convex,
\item $\phi_t$ leaves $\bdry \Sigma$ fixed,
\item $\phi_1(\Gamma_\Sigma)=\Gamma_{\phi_1(\Sigma)}$,
\item $\phi_1(C)$ is Legendrian.
\ee
\end{thm}
Therefore, in particular, a nonisolating collection $C$ can be realized by a Legendrian collection
$C'$ with the same number of geometric intersections.    A corollary of this theorem, observed by
Kanda, is the following:

\begin{cor}[Kanda] A closed curve $C$ on $\Sigma$ can be realized as a Legendrian curve (in the
sense of Theorem \ref{lerp}) if $C\pitchfork \Gamma_\Sigma$ and $C\cap \Gamma_\Sigma\not =
\emptyset$.
\end{cor}

Observe that  if $C$ is a Legendrian curve on a convex surface $\Sigma$, then its twisting
number $t(C,Fr_\Sigma)={1\over 2}\#(C\cap \Gamma_\Sigma)$, where $\#(C\cap \Gamma_\Sigma)$
is the geometric intersection number (signs ignored).

\begin{proof}
By Giroux's Flexibility Theorem, it suffices to find a
characteristic foliation ${\cal F}$ on $\Sigma$ with (an isotopic copy of) $C$ which is represented
by Legendrian curves and arcs.  We remark here that these Legendrian curves and arcs constructed
will always pass through singular points of ${\cal F}$.
Consider a component $\Sigma_0$ of  $\Sigma\backslash (\Gamma_\Sigma \cup
C)$ -- let us assume $\Sigma_0\subset R_+$, so all the elliptic singular points are
sources.
Denote $\bdry \Sigma_0 = \gamma^-- \gamma^+$, where $\gamma^-$ consists of
closed curves $\gamma$ which intersect $\Gamma_\Sigma$, and $\gamma^+$ consists of
closed curves $\gamma\subset C$.
This means that for $\gamma\subset \gamma^-$,
either $\gamma\subset \Gamma_\Sigma$ or $\gamma=
\delta_1\cup \delta_2\cup\cdots \cup \delta_{2k}$, where
$\delta_{2i-1}$, $i=1,\cdots,k$, are subarcs of $C$, $\delta_{2i}$, $i=1,\cdots, k$,
are subarcs of $\Gamma_\Sigma$, and the endpoint of $\delta_j$ is the initial point of $\delta_{j+1}$.
Since $C$ is nonisolating, $\gamma^-$ is nonempty.
What the $\gamma^-$ provide are `escape routes' for the flows whose
sources are $\gamma^+$ or the singular set of $\Sigma_0$  --
in other words, the flow would be exiting  along
$\Gamma_\Sigma$.

Construct ${\cal F}$ so that (1) the subarcs of $\gamma^-$ coming from $C$
are now Legendrian, with a single positive half-hyperbolic point in the interior of the arc,
(2) the curves of
$\bdry \Sigma_0$ contained in $C$ are Legendrian curves, with one positive half-elliptic point
and one positive half-hyperbolic point.   If $\gamma\subset \gamma^-$ intersects $C$,
then we give a neighborhood $\gamma\times I$ a characteristic foliation as in Figure \ref{fig2}.
\begin{figure}
	{\epsfysize=1.5in\centerline{\epsfbox{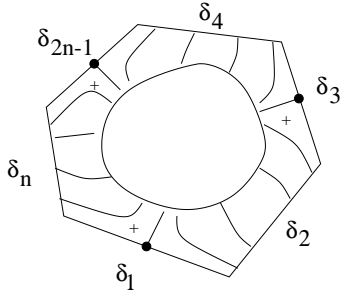}}}
	\caption{Characteristic foliation on $\gamma\times I$}
	\label{fig2}
\end{figure}
After filling in this collar, we may assume that ${\cal F}$ is transverse to and flows out
of $\gamma^-$.   If $\gamma^+$ is empty, then we introduce a positive elliptic singular
point on the interior of $\Sigma_0$, and let $\gamma^+$ be a small closed loop around
the singular point, transverse to the flow.  At any rate, we may assume the flow enters through
$\gamma^+$ and exits through $\gamma^-$ -- by filling in appropriate positive hyperbolic
points we may extend ${\cal F}$ to all of $\Sigma_0$.
\end{proof}

\subsection{Bypasses}

\begin{figure}
	{\epsfysize=1.5in\centerline{\epsfbox{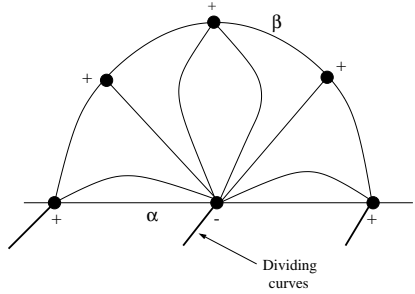}}}
	\caption{A bypass}
	\label{fig3}
\end{figure}

Let $\Sigma$ be a convex surface.
A {\it bypass} is a half-disk $D$ with $\bdry D=\alpha\cup \beta$ Legendrian arcs
such that, for one orientation of $D$, $\alpha\cap \beta$
are both positive elliptic, the other singular point along $\alpha$ is negative elliptic, and all the
singular points along $\beta$ are positive  and alternate between elliptic and hyperbolic.
We assume that all the singular points along $\alpha$ lie on $\Gamma_\Sigma$.
See Figure \ref{fig3}.
Slightly informally, we have twisting numbers $t(\alpha)=-1$ and $t(\beta)=0$;  in other words,
going around the bypass increases the twisting number.  Since it is easy to decrease the
twisting number but not always possible to increase the twisting number, the existence of a bypass
is not a local condition and we do not get bypasses `for free.'

\begin{figure}[ht]
	{\epsfysize=1.5in\centerline{\epsfbox{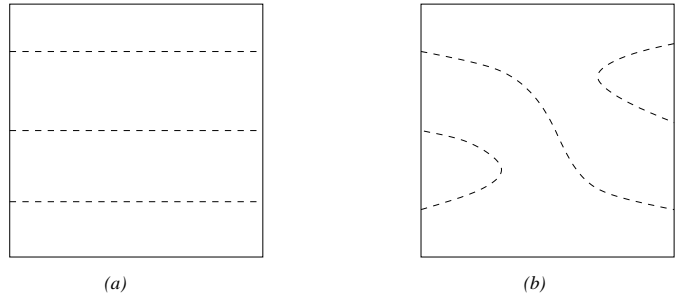}}}
	\caption{(a) Dividing curves (dashed lines) on $\Sigma$, (b) Dividing curves on $\Sigma'$}
	\label{bypassmove}
\end{figure}

We obtain the following result regarding the change in the dividing set when a bypass is attached:

\begin{thm}[\cite{H1}]\label{thm:bypass}
	Let $\Sigma$ be a convex surface and $D$ a bypass on $\Sigma.$
	Then we can find a neighborhood $N$ of $\Sigma\cup D$ with convex
	$\bdry N=\Sigma-\Sigma'$, and $\Gamma_{\Sigma'}$ is  related to
	$\Gamma_\Sigma$ as shown in Figure \ref{bypassmove}.
\end{thm}

The following lemma was observed by W. Kazez:

\begin{lemma}
Let $R_+$ and $R_-$ be the positive and negative regions of $\Sigma$, and
$R_+'$ and $R_-'$ be the positive and negative regions of $\Sigma'$, obtained
from a bypass attachment onto $\Sigma$.  Then
$$\chi(R_+)-\chi(R_-)=\chi(R_+')-\chi(R_-').$$
\end{lemma}

The proof is simple.  Note that $<e(\xi),\Sigma>=\chi(R_+)-\chi(R_-)$, where
$e(\xi)$ is the Euler class of $\xi$.

The consequence of this bypass move is that many questions in contact topology take on
a much more combinatorial appearance -- they can be rephrased into (often nonstandard)
questions about curves on surfaces.  One of the key features is that contact topology is
intimately connected with {\it positive Dehn twists} and the mapping class group
of an oriented surface, as we will see after calculating
some examples of bypass attachments.  Let us assume that the bypass
attachments take place inside a tight contact manifold.

\begin{figure}
	{\epsfysize=2.5in\centerline{\epsfbox{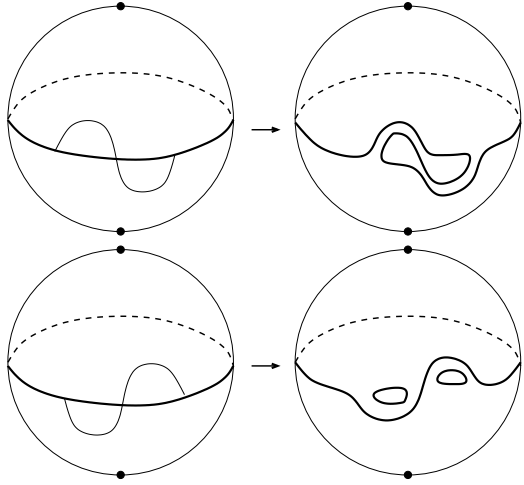}}}
	\caption{Bypasses on $S^2$}
	\label{fig4}
\end{figure}

\vskip.12in
\noindent
{\bf Example:} $\Sigma=S^2$.  There are exactly two possibilities, up to isotopy.
See Figure \ref{fig4}.
The dividing curves are represented by thick lines, whereas the Legendrian curve along which
the bypass is attached is a thin line intersecting $\Gamma_\Sigma$ at three points.
Notice that one of the bypasses preserves $\#\Gamma_\Sigma$, whereas its mirror
image changes $\#\Gamma_\Sigma$ from $1$ to $3$, hence is a disallowed move in
a tight structure.

\begin{figure}
	{\epsfysize=1in\centerline{\epsfbox{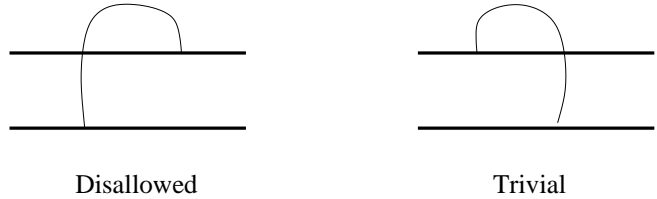}}}
	\caption{Trivial bypass attachment}
	\label{fig5}
\end{figure}

\vskip.12in
In a similar vein,
Figure \ref{fig5} depicts the {\it trivial} bypass attachment, i.e., one that does not change the
dividing curve configuration, together with its evil twin, the disallowed move.   A detailed discussion of
the trivial bypass will be given in Section \ref{Section:trivial}.

\begin{figure}
	{\epsfysize=5in\centerline{\epsfbox{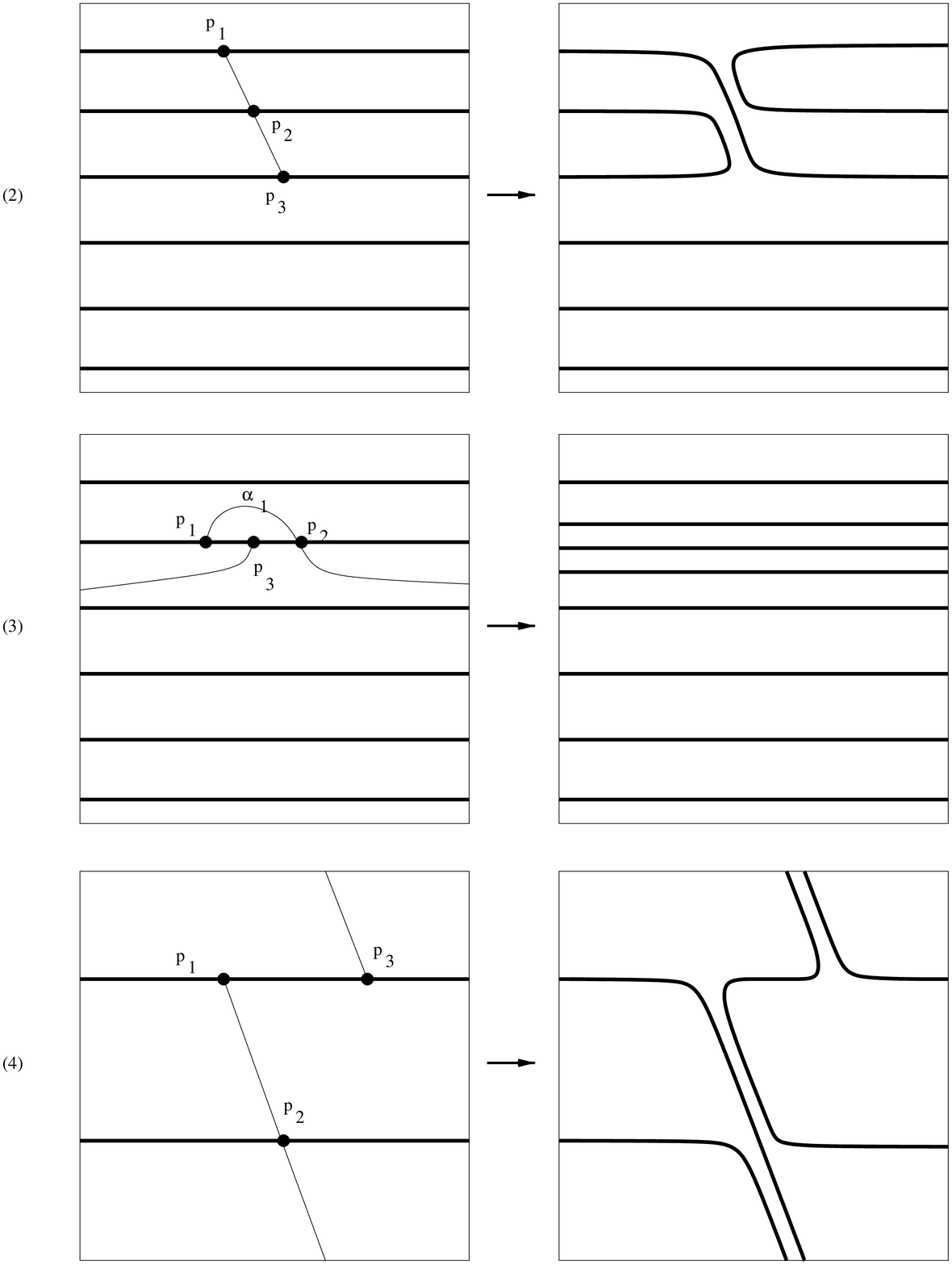}}}
	\caption{Bypass attachments on $T^2$}
	\label{fig6}
\end{figure}

\vskip.12in
\noindent
{\bf Example:}  $\Sigma=T^2$.   The following are the tight possibilities:
\be
\item The bypass attachment is {\it trivial}.
\item $\#\Gamma_\Sigma=2n$ decreases by $2$ (provided $2n>2$).
\item $\#\Gamma_\Sigma=2n$ increases by $2$.
\item The new dividing curve configuration is obtained from performing a (positive)
Dehn twist to the old
configuration.
\ee
The examples are presented in Figure \ref{fig6}.
Let $D$ be the bypass, and $\alpha=D\cap \Sigma$ be the Legendrian arc of attachment.
Let $\{p_1,p_2,p_3\}=\alpha\cap \Gamma_\Sigma$, arranged in order along $\alpha$, and
$\gamma_i$, $i=1,2,3$, be the component of $\Gamma_\Sigma$ containing $p_i$.  If
$\gamma_1=\gamma_2$, then
the subarc $\alpha_s\subset \alpha$ from $p_1$ to $p_2$, together with a subarc
$\gamma_s\subset \gamma_1=\gamma_2$ from $p_1$ to $p_2$ bound a half-disk.  If $p_3\not \in
\gamma_s$, then the bypass attachment is either trivial -- case (1) -- or disallowed.
If $p_3\in \gamma_s$, then the allowed bypass attachment increases $\#\Gamma_\Sigma$
-- this is case (3).  Therefore, we may assume $\gamma_1\not=\gamma_2$ and
$\gamma_2\not=\gamma_3$.  If $\#\Gamma_\Sigma>2$,
then the $p_i$ lie on three distinct dividing curves, and we have case (2).  Finally, we have
$\#\Gamma_\Sigma=2$, $\gamma_1=\gamma_3\not=\gamma_2$.
We now exhibit the curve that
$\Sigma$ is Dehn twisted about, to obtain $\Gamma_{\Sigma'}$ from $\Gamma_\Sigma$.
Since $\gamma_1=\gamma_3$, $\gamma_1\backslash \{p_1,
p_3\}$ consists of two arcs.  Take the arc $\gamma_s$, starting at $p_1$ in the direction
$\gamma_s'$ which forms an oriented basis $\{\gamma_s',\alpha'\}$, where $\alpha'$ is
a tangent vector to $\alpha$ at $p_1$ which points outward.
Now let the curve $C$  for the Dehn twist be $\alpha\cup \gamma_s$.

\begin{prop}
Let $\Sigma$ be a closed convex surface in a tight manifold.  A bypass attachment to $\Sigma$
give rise to $\Sigma'$ which satisfies one of the following:
\be
\item $\Gamma_\Sigma=\Gamma_{\Sigma'}$, i.e., a trivial attachment.
\item $\#\Gamma_\Sigma= \#\Gamma_{\Sigma'}+2$.
\item $\#\Gamma_\Sigma= \#\Gamma_{\Sigma'}-2$.
\item $\Gamma_{\Sigma'}$ is obtained from $\Gamma_\Sigma$ via a positive Dehn twist.
\item `Mystery move', described below.
\ee
\end{prop}

\begin{proof}    The argument is almost identical to the $T^2$ case.  Using the notation from above,
let  $p_i$, $i=1,2,3$, be the three intersections of $\alpha\cap \Gamma_\Sigma$ and
$\gamma_i$ be components of $\Gamma_\Sigma$ containing $p_i$.  As above, if
$\gamma_1\not=\gamma_2$ and $\gamma_2\not=\gamma_3$, then we have a
decrease (2) or a positive Dehn twist (4).  Suppose $\gamma_1=\gamma_2$.
The chief difference between $T^2$ and
higher genus $\Sigma$ is that the arc $\alpha_s$ from $p_1$ to $p_2$ does not always
bound
a half-disk, together with an arc $\gamma_s\subset \gamma_1$ -- if it does, then we have
(1) or (3).
We have two remaining cases: $\gamma_1=\gamma_2\not=\gamma_3$, which we call
the `mystery move' (5),
and $\gamma_1=\gamma_2=\gamma_3$, which gives either (2) or (4).    See Figure
\ref{fig9} for these possibilities.
\begin{figure}
	{\epsfysize=3.5in\centerline{\epsfbox{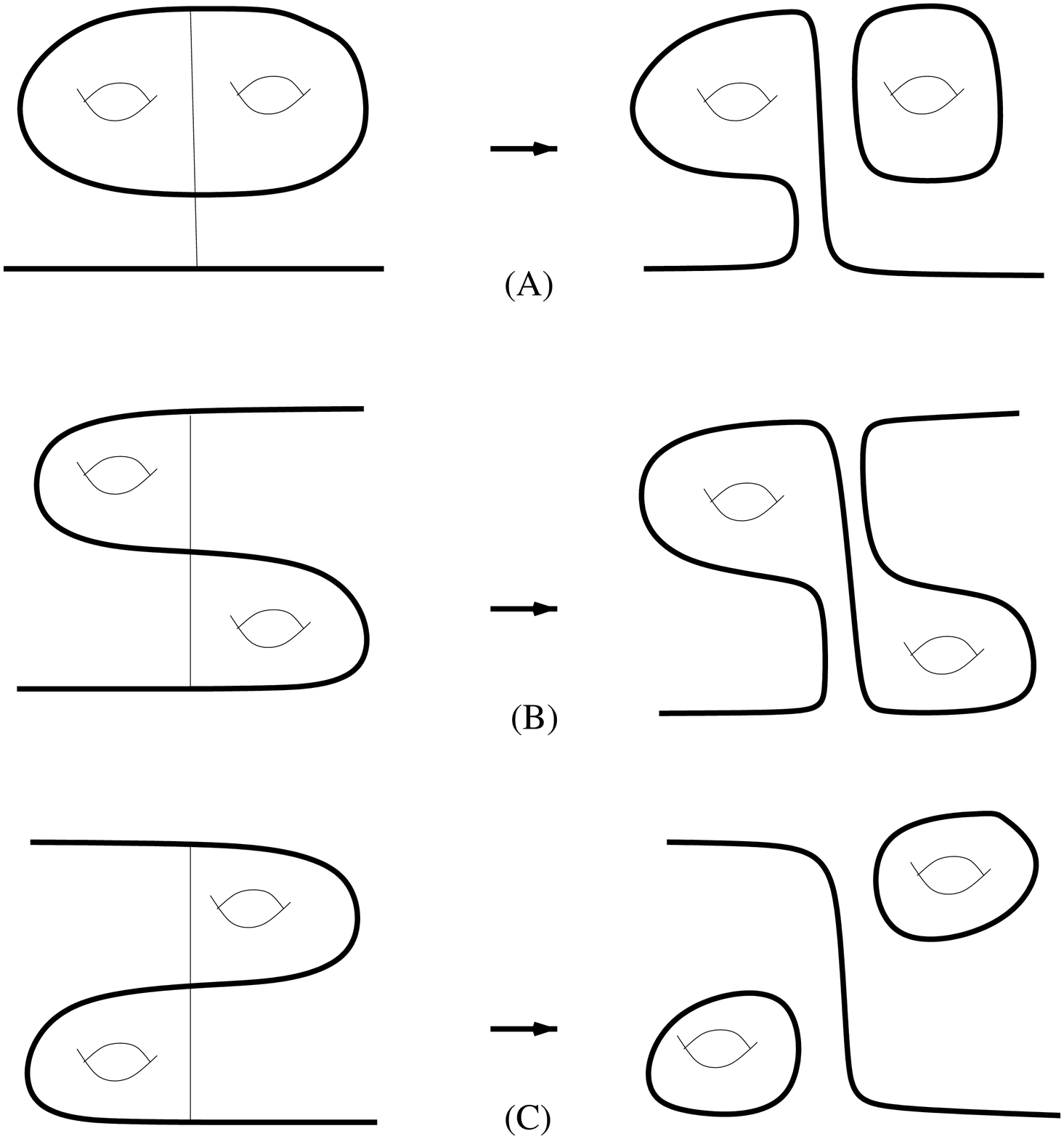}}}
	\caption{Bypasses on $\Sigma$: (A) is the mystery move, (B) is given by a Dehn twist,
and (C) increases $\#\Gamma_\Sigma$ by 2}
	\label{fig9}
\end{figure}
\end{proof}

Positive Dehn twists have also appeared in the study of symplectic Lefschetz fibrations, so
it is not as surprising that they appear in contact topology.  Understanding the relationship among
the following:
\vskip.12in
$$ \begin{array}{c}
\mbox{Contact topology}\\
\updownarrow\\
\{\mbox{Positive Dehn twists}\}\subset \mbox{Mapping class group}\\
\updownarrow\\
\mbox{Lefschetz fibrations}
\end{array}$$
\vskip.12in
\noindent
is an interesting open-ended question.

\vskip.15in
Let us now examine several examples where $\Sigma$ is compact convex with
Legendrian boundary.
\vskip.15in
\noindent
{\bf Example:} $\Sigma=D^2$.  Assume $tb(\bdry \Sigma)=-n<0$.  Since there cannot
exist closed dividing curves on $\Sigma$, all the dividing curves are arcs which connect
two points on the boundary.  There are $2n$ points which need to be `hooked up' via
the dividing curves, so the number of possible configurations is finite -- largely due
to the nonexistence of homotopically nontrivial closed curves.  As before,
let $\alpha$ be the Legendrian arc of attachment, and $p_i$ be its intersection with
$\gamma_i$, $i=1,2,3$.
There are two possibilities for bypass attachments:
\be
\item The bypass attachment is trivial $\Leftrightarrow$
$\gamma_1=\gamma_2$ or $\gamma_2=\gamma_3$.
\item $\gamma_1,\gamma_2, \gamma_3$ are distinct.
\ee

\begin{figure}
	{\epsfysize=1.5in\centerline{\epsfbox{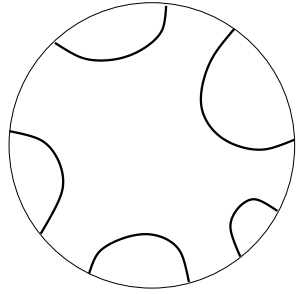}}}
	\caption{No nontrivial bypasses}
	\label{fig7}
\end{figure}

There is one configuration where (2) cannot happen.  This is when each dividing curve, together
with an arc on $\bdry \Sigma$,
bounds a half-disk component of $\Sigma\backslash \Gamma_\Sigma$.
Refer to Figure \ref{fig7}.

\vskip.1in
Given a convex surface $\Sigma$ with Legendrian boundary,  a dividing curve $\delta$ is called
{\it boundary-parallel} ($\bdry$-parallel) if $\delta$ cuts off a half-disk which has no other intersections
with $\Gamma_\Sigma$.   A dividing set $\Gamma_\Sigma$ is called {\it boundary-parallel} if
all its dividing curves are $\bdry$-parallel arcs.

\begin{figure}
	{\epsfysize=1.5in\centerline{\epsfbox{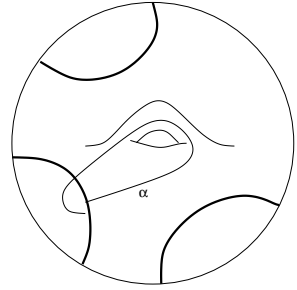}}}
	\caption{Increasing $\#\Gamma_\Sigma$}
	\label{fig8}
\end{figure}

\vskip.15in
\noindent
{\bf Example:} $\Sigma$ is any compact convex surface with Legendrian boundary.
If $\Gamma_\Sigma$ is $\bdry$-parallel, then the only
nontrivial bypass attachment is one which increases $\#\Gamma_\Sigma$.  See Figure \ref{fig8}.

\subsection{Discussion of trivial bypasses}   \label{Section:trivial}
In this section we prove two lemmas which guarantee the existence and triviality of the
so-called {\it trivial} bypasses in the previous section.

\begin{lemma}[Existence] \label{existence}
Let $\Sigma$ be a convex surface which is closed or compact with Legendrian
boundary, and assume the ambient contact manifold $(M,\xi)$ is tight.
Let $\delta\subset \Sigma$ be any Legendrian arc, drawn as in the right-hand diagram of Figure \ref{fig5}.
Then there exists a trivial bypass attached along $\delta$, which lies entirely in an $I$-invariant neighborhood
of $\Sigma$.
\end{lemma}

\begin{proof}
The key ingredient is the Legendrian realization principle.  We apply it to $\Sigma$, fixing
$\Gamma_\Sigma$ and $\delta$, so that we may assume that there exists a convex disk
$D\subset \Sigma$ with collared Legendrian boundary, $tb(\bdry D)=-2$, $\#\Gamma_D=2$,
and $\delta\subset D$.  All the operations which follow will now take place inside an $I$-invariant
neighborhood $D\times I$, where $D=D\times\{0\}$. Let $p,q$ be two points on
the same component of $\Gamma_D$.  On $D\times \{0\}$, take a Legendrian arc $\delta_0$
with endpoints on $\{p,q\}\times\{0\}$
which extends $\delta$ as in the left-hand diagram of Figure \ref{fig16}.
\begin{figure}
	{\epsfysize=1.5in\centerline{\epsfbox{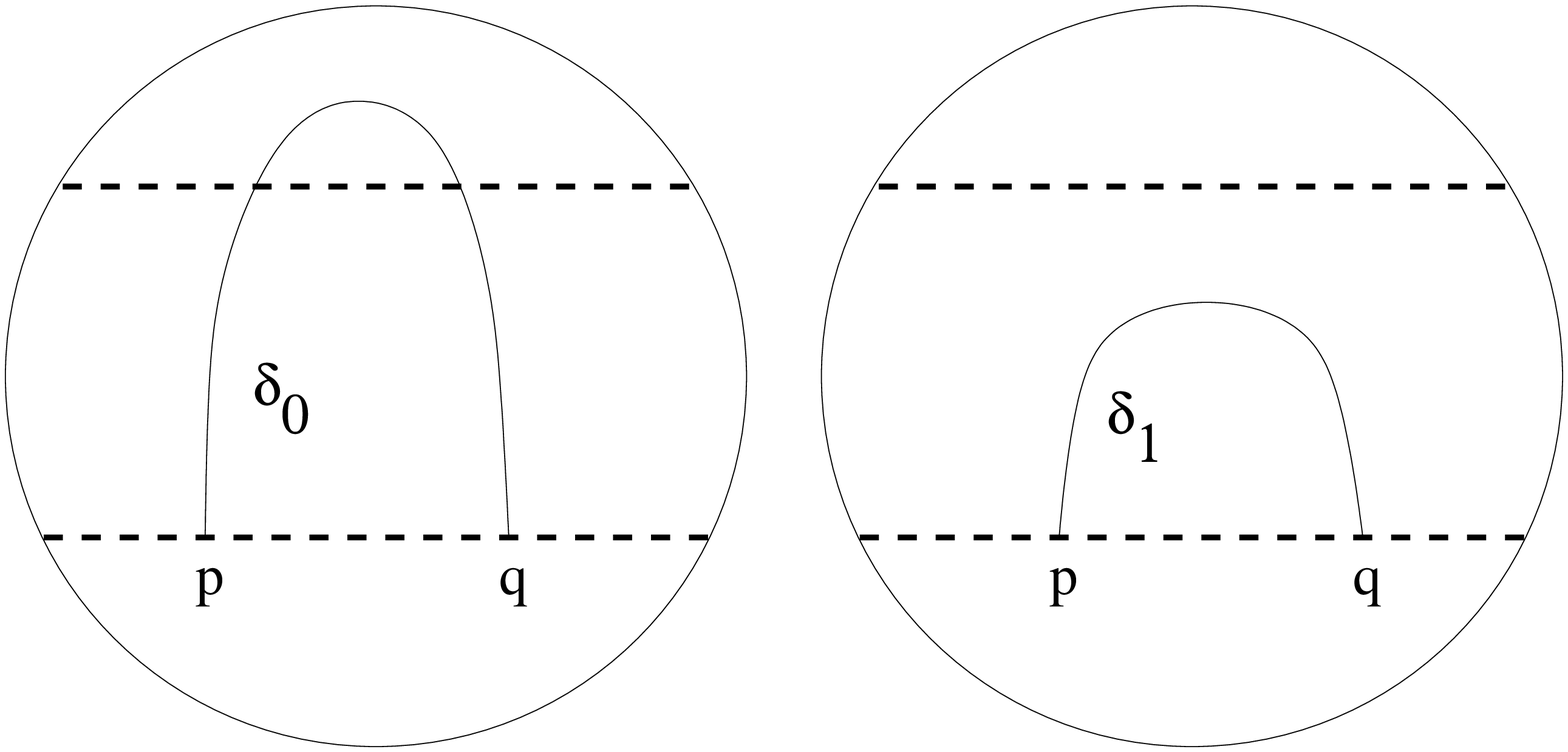}}}
	\caption{}
	\label{fig16}
\end{figure}
On $D\times \{1\}$, take another Legendrian arc $\delta_1$ depicted in the right-hand diagram of
Figure \ref{fig16}, with endpoints $\{p,q\}\times\{1\}$.  Note that
both $D\times \{0\}$ and $D\times \{1\}$ may be slightly modified using Giroux's Flexibility Theorem or the Legendrian
realization principle to realize $\delta_0$, $\delta_1$ as Legendrian arcs.
Now, form the closed Legendrian curve $\gamma=\delta_0\cup \delta_1\cup  (\{p,q\}\times I)$.
Let $D'$ be a convex disk with Legendrian boundary $\gamma$.  Since $tb(\gamma)=-2$, there are two
possibilities for bypasses along $\delta_0$, namely the trivial bypass and the disallowed bypass.  The
attachment must be trivial since $\xi$ is tight.
\end{proof}

The next lemma proves that a trivial bypass is, indeed, trivial --- a trivial bypass may be attached to
a convex surface without effect.

\begin{lemma} [Triviality] \label{triviality}
Let $\Sigma$ be a convex surface which is closed or is compact with Legendrian
boundary.  If a trivial bypass $D'$ is attached along $\delta\subset \Sigma$ as in the right-hand side of
Figure \ref{fig16}, then a neighborhood $N(\Sigma\cup_\delta D')$ is isotopic to the standard $I$-invariant
neighborhood of $\Sigma$.
\end{lemma}

\begin{proof}
As in the Existence Lemma, isolate $\delta$ inside a convex disk $D\subset \Sigma$ with collared Legendrian
boundary and $tb(\bdry D)=-2$, after possible perturbation of $\Sigma$, fixing $\Gamma_\Sigma$
and $\delta$.    Here, we modify the characteristic foliation on $D$ rel $\Gamma_D\cup \delta$, so that
it matches that of the local model in the Existence Lemma.  Now, a thickened neighborhood
$N(D\cup_\delta D')$ of $D\cup_\delta D'$ is tight since a model exists.
Moreover, it is $I$-invariant due to the uniqueness of
the tight contact structure on the 3-ball with fixed boundary (Theorem \ref{thm:unique}).
Finally, we conclude that $N(\Sigma\cup_\delta D')$ is $I$-invariant.
\end{proof}

\subsection{Convex decompositions}
One of the effective ways of decomposing a tight contact 3-manifold is to successively cut along
convex surfaces. The possibility of performing a Haken decomposition was already present
in Kanda's seminal paper \cite{K97}, where the proof of the classification of the 3-torus $T^3$
depended on a Haken decomposition.   Another hint of a deeper connection appeared in
Torisu's paper \cite{To}, where the analogy between sutures for foliations and dividing curves for
contact structures are made.    The author learned about the possibility of convex decompositions
from John Etnyre.  In this section, we present a glimpse of this theory.
These ideas will be developed in a joint paper with Will Kazez and Gordana Mati\'c.
In particular, we will reprove a result of Eliashberg and Thurston on the
existence of a universally tight contact structure on any 3-manifold $M$ with
$H_2(M,\bdry M)\not=0$, by directly using convex decompositions and the gluing
theorem in the next section.

If the surface of the cutting has boundary, then we ask that the boundary
be made Legendrian, so that the Haken decomposition occurs along a convex surface with Legendrian
boundary.   More precisely, let $(M,\xi)$ be a tight contact manifold with convex
boundary $\bdry M=\Sigma$
and we cut along $S$ with boundary on $\Sigma$.  Let $\gamma$ be a boundary component
of $\bdry S\subset \Sigma$.   We choose $\gamma$ in its isotopy class so that the geometric
intersection
$\#(\gamma\cap \Gamma_\Sigma)$ is minimized, provided this number is $\geq 2$.  If the minimum
geometric intersection is $0$, then choose $\gamma$ so $\#(\gamma\cap \Gamma_\Sigma)=2$  --
we artificially force this because cutting along Legendrian curves with twisting number $0$ is not
as easy to control.
(Note $\Gamma_\Sigma$ is never empty.)  Then use the Legendrian realization principle to make
$\bdry S$ Legendrian and $S$ convex.  Once we have prepared $S$ as above, we can perform the
splitting.  In order to ensure the convexity of the resulting surface, we apply edge-rounding,
described in the next paragraph.

Let $\Sigma_1$ and $\Sigma_2$
be compact convex surfaces with Legendrian boundary, which
intersect transversely along a common boundary Legendrian curve $L$.
The neighborhood of
the common boundary Legendrian is locally isomorphic to the neighborhood $\{x^2
+ y^2  \le \varepsilon\}$ of $M = {\bf R}^2 \times ({\bf R}/{\bf Z})$ with
coordinates $(x, y, z)$ and con\-tact 1-form $\alpha = \sin(2\pi n z) dx +
\cos(2\pi n z) dy$, for some $n \in {\bf Z}^+$.  Let $A_i\subset \Sigma_i$, $i=1,2$,
be an annular collar of the boundary component $L$.
We may choose our local model so that $A_1 = \{x = 0, 0 \le
y \le \varepsilon\}$
and $A_2 = \{y = 0, 0 \le x \le \varepsilon\}$ (or the same with $A_1$ and $A_2$  switched).
Assuming the former, if we join $\Sigma_1$
and $\Sigma_2$ along
$x = y = 0$ and round the common edge, the resulting surface is convex, and the
dividing curve
$z = \frac{k}{2n}$ on $\Sigma_1$ will connect to the dividing curve $z = \frac{k}
{2n} - \frac{1}
{4n}$ on $\Sigma_2$, where $k = 0, \cdots, 2n - 1$.

Hence we will be obtaining a decomposition
$$ M=M_0\sa M_1 \sa \cdots \sa M_n=\cup B^3.$$
Note that the end product of the decomposition should be a union of 3-balls.
If the contact structure is tight, then $\bdry B^3$ should all have exactly one dividing curve.

Let us now recall the following fundamental theorem of Eliashberg \cite{E92}:
\begin{thm}[Eliashberg]  \label{thm:unique} Assume there exists a contact structure $\xi$ on a neighborhood of
$\bdry B^3$ which makes $\bdry B^3$ convex with $\#\Gamma_{\bdry B^3}=1$.
Then there exists a unique extension of $\xi$ to a tight contact structure on $B^3$
up to an isotopy which fixes the boundary.
\end{thm}

The basic building blocks are $B^3$, each with a unique tight contact structure.
The main goal of this paper is to explore how to glue the manifold back along the
convex decomposition to obtain a tight contact structure.

We briefly mention the relationship between convex decompositions and
sutured manifold decompositions of Gabai \cite{Ga}.
A {\it sutured manifold} $(M,\gamma)$ is a compact oriented 3-manifold together
with $\gamma\subset \bdry M$ consisting of pairwise disjoint annuli $A(\gamma)$
and tori $T(\gamma)$.  A suture $\gamma$ on $\bdry M$ {\it divides} $\bdry M\backslash \gamma$
into positive and negative regions $R_+$ and $R_-$, and every annular suture bounds a component of
$R_+$ on one side and a component of $R_-$ on the other.
We now define a  {\it sutured manifold decomposition} $(M,\gamma) \stackrel{S}{\sa}
(M',\gamma')$.
First assume $S$ satisfies:
\be
\item Every component of $\bdry S$ is a homotopically nontrivial curve in $\bdry M$.
\item $\bdry S \pitchfork \gamma$.
\item Every disk component of $S$ intersects $\gamma$.
\ee
Let $M'=M\backslash S$,
$R'_+=(R_+\backslash S)\cup S_+$ and $R'_-=(R_-\backslash S)
\cup S_-$, where
$S_\pm$ are the parallel copies of $S$ on the cut-open manifold $M'$.
Define $\gamma'$ to consist of the common boundary of $R'_+$ and $R_-'$, together with
the remaining $T^2$ components -- the  $T^2$ components which are cut become annuli, and
each annulus is squashed to its core curve.

The following correspondence shows that a sutured manifold decomposition is a special
case of the convex decomposition:
\be
\item A annular suture can be viewed as a dividing set, if we squash the annulus to its core curve.
A sutured manifold decomposition along a $T^2$ component can be interpreted as follows:
Right before we cut along $S$ with $S\cap T^2\not=\emptyset$, we substitute the
$T^2$ suture by a pair of parallel homotopically nontrivial dividing curves, each of  which has
algebraic intersection $1$ with
each component of $S\cap T^2$.
\item A component $\Sigma\subset \bdry M$ may not have a suture at all, whereas a dividing
set must not be empty.  We remedy this by placing a pair of parallel homotopically
nontrivial dividing curves on $\Sigma$  with nontrivial geometric intersection with each component
of $\bdry S$,
before cutting.
\item Let  $S$ be a convex surface -- realize the boundary as a Legendrian curve with
twisting number $\leq -2$ -- and choose $\Gamma_S$ so that every dividing curve is an
arc which is $\bdry$-parallel.
\item When $M$ is cut along $S$ and rounded, all the dividing curves of the decomposed
manifold $M'$, except perhaps for
the $T^2$ components and components $\Sigma\subset M$ without sutures, correspond
to sutures  $\gamma'$ on $M'$.
\ee

\section{Gluing theorem}

\subsection{Statements}

Let $M$ be an oriented, compact 3-manifold with nonempty boundary.  Fix a dividing set
$\Gamma_{\bdry M}$ and some singular foliation $\mathcal{F}$ which is {\it adapted to}  $\Gamma_{\bdry M}$.
Define $\mathcal{T}(M,\mathcal{F})$ to be the set of isotopy classes of tight contact structures on
$M$ which restrict to $\bdry M$ to a characteristic foliation $\mathcal{F}$.  Using Giroux's Flexibility Theorem,
it is easy to see that there is a natural bijection
$$\mathcal{T}(M,\mathcal{F})\stackrel{\sim}{\rightarrow} \mathcal{T}(M,\mathcal{F'}),$$
for every pair of singular foliations $\mathcal{F}$ and $\mathcal{F'}$ adapted to $\Gamma_{\bdry M}$.
(For more details, see \cite{H1}.)  Therefore, we write $\mathcal{T}(M,\Gamma_{\bdry M})$ to mean
any of the $\mathcal{T}(M,\mathcal{F})$.

\subsubsection{Handlebody Case}    \label{handlebody}
Let $M$ be a genus $g$ handlebody, and $\xi$ a contact structure on $M$ so that
$\Sigma=\bdry M$ is convex.  Assume
$D_1,\cdots,D_g$ are compressing disks which give $M\backslash (D_1\cup \cdots \cup D_g)=B^3$.
We may choose $\bdry D_i$ to be a Legendrian curve on $\Sigma$, using the
Legendrian realization principle (so we are choosing a particular characteristic foliation on
$\Sigma$ consistent with $\Gamma_\Sigma$).  We then ask all the $D_i$ to be
convex with Legendrian boundary. If $\xi$ is tight, then we require $tb(\bdry D_i)=n_i<0$.

The characteristic foliation (or, better, the dividing curves) on the 2-skeleton,
consisting of convex surfaces $\Sigma$ and $D_1,\cdots,D_g$, uniquely determines the
contact structure on the rest of $M$, if $\xi$ is to be tight.   This follows from Eliashberg's
uniqueness theorem (Theorem \ref{thm:unique}), applied to the cut-open manifold
$M\backslash (D_1\cup \cdots \cup D_g)$.
Hence we encode all the information for $M$ via dividing curves on $\Sigma$ and
the $D_i$.

Now, let  $\mathcal{T}^*(M,\Gamma_{\bdry M})=\mathcal{T}(M,\Gamma_{\bdry M})\cup \{*\}$, where $*$ is a single point
corresponding to all the overtwisted contact structures.  Also define the
{\it configuration space} ${\mathcal C}$ to be the set of $C=(\Gamma_1,\cdots,
\Gamma_g)$, where $\Gamma_i$ is a possible dividing set for a convex disk $D_i$
with $tb(\bdry D_i)=n_i$, subject to the condition that $D_i$ has a tight neighborhood.
By Giroux's criterion, this is equivalent to saying that $\Gamma_i$ has no
closed curves on $D_i$, since any closed curve must necessarily bound a disk. Then there is a map:
$$\Psi: \mathcal{C}\rightarrow \mathcal{T}^*(M,\Gamma_{\bdry M})$$
which sends a configuration $C$ to its corresponding tight contact structure $\Psi(C)$
if the glued-up contact structure is tight, and to $*$ if the glued-up contact structure is
overtwisted.      $\Psi$ is surjective, but not necessarily injective, since it is possible that a tight contact structure
$\xi$ could have arisen from multiple configurations.

To remedy this situation, we introduce a directed graph $G=({\mathcal C},{\mathcal T})$, where the configuration
space $\mathcal{C}$ is the set of vertices and
${\mathcal T}\subset {\mathcal C}\times
{\mathcal C}$ is the set of edges or
{\it allowable state transitions}, which we define in the next few sentences.  We often write
$C\sc C'$ for $(C,C')\in \cal T$.

Write $B=M\backslash(D_1\cup\cdots \cup D_g)$.
A configuration $C=(\Gamma_1,\cdots,\Gamma_g)\in {\mathcal C}$
gives rise to $\Gamma_{\bdry B}$, after rounding
edges.  We call $C$  {\it potentially allowable}, if $\#\Gamma_{\bdry B}=1$.
Note that $\#\Gamma_{\bdry B}=1$
is equivalent to saying that the contact structure on $\bdry B$ can
be extended uniquely to a tight contact structure on $B$.   We say $C$ is `potentially'
allowable because we do not immediately spot an overtwisted disk.

We say a {\it state transition} is  {\it allowable} and write $C\sc C'$ if
\begin{enumerate}
\item $C$ is potentially allowable.
\item $C'$ can be obtained from $C$ via a single {\it nontrivial} bypass attachment.
\item Attaching the bypass onto $\bdry B$ from the interior of $B$
does not change $\#\Gamma_{\bdry B}$.
\end{enumerate}
It is easy to verify that $C\sc C'$ implies $C'\sc C$, except when $C'$ is already not
potentially allowable.

Now, $C\in \cal C$ is {\it allowable} if every $C'\in \cal C$ in the same `connected component' of
$G$ is {\it potentially allowable}.  In other words, every $C'$ that can be reached via a
sequence of
allowable state transitions, starting at $C$, must not be `obviously overtwisted' when $M$ is
cut open along $\cup_i D_i$ with configuration $C'$.   Denote the
set of allowable $C$ by $\mathcal{C}_0\subset \mathcal{C}$.  On $\mathcal{C}_0$ the graph is reflexive,
and we write $\pi_0(\mathcal{C}_0)$ to mean the connected components of $\mathcal{C}_0$.

\begin{thm}[Gluing/Classification] \label{gluing}
$\Psi$ restricts to a surjective map $\mathcal{C}_0\rightarrow \mathcal{T}(M,\Gamma_{\bdry M})$, which
in turn factors through $$ \pi_0(\mathcal{C}_0)\stackrel {\sim}{\rightarrow} \mathcal{T}(M,\Gamma_{\bdry M}).$$
\end{thm}

We have the following corollary:

\begin{cor}
Let $[C]\in \pi_0(\mathcal{C}_0)$ be the connected component containing $C$.  If $[C]$ contains only
one configuration (i.e., there are no state transitions from $C$), then the corresponding
contact structure is universally tight.
\end{cor}

It is also relatively easy to show the following corollary using state transitions:

\begin{cor}
A tight contact structure on a solid torus with convex boundary is universally tight if and only if
the corresponding $[C]$ consists of exactly one element.
\end{cor}

We end this section with the following question:

\begin{q}
Characterize, in terms of state transition data, what it means for a contact structure on a
genus $g$ handlebody to be universally tight.
\end{q}

\subsubsection{General case}
The Gluing/Classification Theorem can be stated in greater generality -- the catch is that the conditions are difficult
to verify in general when the cutting surfaces are not disks.
Let $M$ be a compact oriented 3-manifold.
If $\bdry M\not= \emptyset$, then we prescribe
a {\it tight} $\Gamma_{\bdry M}$, i.e., a dividing set which comes from a
tight contact structure on a neighborhood of $\bdry M$.
Cut $M$ along an incompressible surface
$N$, possibly with Legendrian boundary to obtain $M\stackrel{N}\sa M'$.
Let us assume that $\bdry N\cap \Gamma_{\bdry M}
\not = \emptyset $ if $\bdry N \not =\emptyset$.
Denote by ${\cal C}$ the {\it configuration space}, consisting of $(\Gamma_N,\xi)$,
where $\Gamma_N$ is a dividing set of $N$,  $\xi$ is an isotopy class of  contact
structures (rel boundary) on $(M',\Gamma')$, and $\Gamma'$ is obtained from $\Gamma_N$.
(At this point we do not assume $\Gamma_N$ is tight or $\xi$ is tight.)
We say $(\Gamma_N,\xi)$ is {\it potentially allowable},
if $\xi$ on $M'$ is tight.   A state transition $C=(\Gamma_N,\xi)\sc
C'=(\Gamma_N',\xi')$ is {\it allowable} if:
\be
\item $C$ is potentially allowable.
\item $C'$ is obtained from $C$ via a single nontrivial bypass attachment along $N$.
\item $\xi'$ is obtained from $\xi$ by peeling an $N\times I$ layer from $M'$, corresponding to
the bypass attachment, (along one copy of $N$),  and reattaching the layer along the other
copy of $N$ on $M'$.
\ee
Also, $C$ is {\it allowable} if every $C'$ in the same connected component is
potentially allowable.    Again, define $\mathcal{C}_0\subset \mathcal{C}$ to be the subset of allowable
$C$, and let $\pi_0(\mathcal{C}_0)$ be the connected components of $\mathcal{C}_0$.

Define ${\cal T}(M,\Gamma_{\bdry M})$  to be the set of isotopy classes of
tight contact structures on $M$
with boundary condition $\Gamma_{\bdry M}$ -- assume the isotopy is rel the boundary if the
boundary is nonempty.    (Also define ${\cal T}(M)$ to be the set of isotopy classes of tight contact structures on
$M$.)
The theorem is then:

\begin{thm}[Gluing/Classification] \label{classification}  Let $M$ be a compact, oriented, irreducible
3-manifold, and
$M\stackrel{N}{\sa} M'$ be a decomposition along an incompressible surface $N\subset M$.
If $\bdry M\not= \emptyset$, then we prescribe
$\Gamma_{\bdry M}$, and, if $\bdry N\not=
\emptyset$, then let $\bdry N$ be Legendrian with $t(\gamma)<0$ for each component
$\gamma\subset \bdry N$.    Then
${\cal T}(M,\Gamma_{\bdry M})$ is in $1-1$ correspondence with $\pi_0(\mathcal{C}_0)$.
\end{thm}

The incompressible surface $N$ does not need to be connected.
The irreducibility of each component of $M\backslash N$ is a useful assumption which facilitates the
proof of the theorem.  However, in the case $M$ has a maximal connect sum decomposition
$M=M_1\# \cdots\# M_n$, where each capped-off manifold $\overline{M_i}=M_i\# B^3$ is
irreducible, we may still apply the same proof with  $N=\cup S^2$ first, and then analyze each
of the irreducible components $\overline{M_i}$.  See Corollary 1 below.

A repeated application using a convex decomposition (until we get a union of $B^3$'s)
gives a complete classification theorem for tight contact structures.
Unfortunately, in every case besides the disk,
where there only finitely many states,
there usually are infinitely many states, which makes the combinatorial problem an infinite one.

The following corollaries represent combinatorially trivial cases of Theorem \ref{classification}.
(Corollary 1 is a consequence of the proof of Theorem \ref{classification}.)
Observe that if $C$ is a potentially allowable configuration, and there
are no state transitions from $C$, i.e., there are no nontrivial bypass attachments, then the corresponding
glued-up contact structure must be tight.   Recall that a tight $\Gamma_{S^2}$ is unique and $\#\Gamma_{S^2}
=1$. For cases 1 and 2 below, we have already discussed that all bypasses which can be attached are trivial.
The proof of the third corollary will appear in a separate paper.

\begin{cor}
\be
\item[] \vskip.1in
\item  (Colin \cite{Co97}, Makar-Limanov \cite{ML})
Consider the connect sum decomposition
$M$ $=$ $M_1\#$ $\cdots$ $\# M_n$, where each capped-off manifold $\overline{M_i}=M_i\# B^3$ is
irreducible.  Then there is a bijection
$$\mathcal{T}(M)\stackrel{\sim}{\rightarrow} \mathcal{T}(\overline{M_1})\times \cdots \times \mathcal{T}(
\overline{M_n}).$$

\item (Essentially due to Colin \cite{Co99})
If we glue a tight contact manifold along disks $D_1$, $D_2$
which are convex with Legendrian boundary, and $\Gamma_{D_i}$, $i=1,2$,
are $\bdry$-parallel, then the glued manifold is tight.
\item (Essentially due to Colin \cite{Co99})
Consider a convex decomposition $(M,\xi)\stackrel{\Sigma}{\sa}(M',\xi|_{M'})$,
where $\xi|_{M'}$ is universally tight, $\Sigma$ is convex with
Legendrian boundary, $\Gamma_{\Sigma}$
is $\bdry$-parallel, and $\Sigma$ $\pi_1$-injects into $M$. Then $(M,\xi)$ is universally tight.
\ee
\end{cor}

\subsection{Proofs of Theorems \ref{gluing} and \ref{classification}}
In this section we first prove Theorem \ref{classification}, and obtain Theorem
\ref{gluing} as a consequence.   We can think of Theorem
\ref{classification} as consisting of two parts: gluing and distinguishing.

\vskip.12in
\noindent
{\bf Idea of Gluing:}  Let $\xi$ be a contact structure on $M$,
and $N$ an incompressible convex surface in $M$.  Consider an overtwisted disk $D\subset M$.
If $D\cap N=\emptyset$, we have a contradiction, because we assume $M\backslash N$ is
tight.  Therefore, $D\cap N\not=\emptyset$.  We carefully control the intersection
$D\cap N$, and show that we may detach $D$ and $N$ after a sequence of bypass moves.

\subsubsection{Controlled intersection of overtwisted disk with convex surface}

\begin{lemma}
Let $D$ be an overtwisted disk in $M$, and $N$ be a convex surface in $M$.
After a possible contact isotopy, we may assume that
$D\pitchfork \Sigma$ and $\bdry D\cap \Sigma\subset \Gamma_\Sigma$.
\end{lemma}

\begin{proof}
Perturb and push $\bdry D$ off the {\it Legendrian skeleton} $K$ of $N$ and
make $D$ transverse to $N$.
The set $K$ consists of the singular points and trajectories connecting between
singular points of the same sign, assuming that
the flow on $N$ has been perturbed into a Morse-Smale one without closed orbits --
$K$ is a deformation retract of $N\backslash \Gamma_N$.   We can then push $\bdry D$
into the dividing set as follows:
Along the trajectory from $K$ to $\Gamma_N$ containg $p\in \bdry D$, the contact structure $\xi$
has the form $\sin y \mbox{ } dz+\cos y \mbox{ } dx$, on $\R^3$ with
coordinates $(x,y,z)$, where $\{-1\leq x\leq 1, -{\pi\over 2}\leq y\leq 0\}$ is the
neighborhood in $N$ of the trajectory $x=z=0$ containing $p$, and
$y=0$ corresponds to the dividing set of $N$.     Then consider the
front projection onto $\R^2$ using the $x$, $z$ coordinates.  Near $p$, we may
assume $\bdry D$ is a straight line $L=\{z=sx\}$ through $(0,0)$ with negative slope $s$.  Modify $L$
fixing endpoints, to $z=f(x)$ with $f'(x)>0$ everywhere except for $x=0$, where the smooth
curve has vertical slope.  Since we may assume that the $\bdry D\cap N$ occur on distinct
trajectories of $N$, we are done after perturbing $D$ again to make it transverse to $N$.
\end{proof}

Let us now examine $D\cap N\subset D$.  This is a union of closed curves (all homotopically
trivial) and arcs connecting between two points on $\bdry D$.  $D\cap N\subset N$ consists of
closed curves and arcs which connect between $\Gamma_N$ --- the closed curves are
all homotopically trivial on $N$ since $N$ is incompressible.

\begin{lemma}
$D\cap N$ can be made Legendrian on $N$, fixing endpoints of arcs,
after possible modification of $D$.
\end{lemma}

\begin{proof} First note that a closed curve $C$ of $D\cap N$ on $N$ cannot be homotopically
nontrivial  -- this is because $C$ bounds a subdisk of $D$, contradicting the
incompressibility of $N$.  Next, assume $C\cap \Gamma_N=\emptyset$
and $C$ bounds a disk $D'$ in $N$.   If $C$ is the innermost intersection
$D\cap N$ on $D'$, then $D$ can be capped off
to reduce the number of components of $D\cap N$.  If there exists an arc of $D\cap N$
on $D'$, then $C\cap \Gamma_N$ could not have been empty.  Therefore, the conditions
of the Legendrian realization principle are satisfied, and $D\cap N$ can be made Legendrian
while keeping the endpoints of arcs fixed.
\end{proof}

\subsubsection{Pushing disks across}
\begin{figure}
	{\epsfysize=1.5in\centerline{\epsfbox{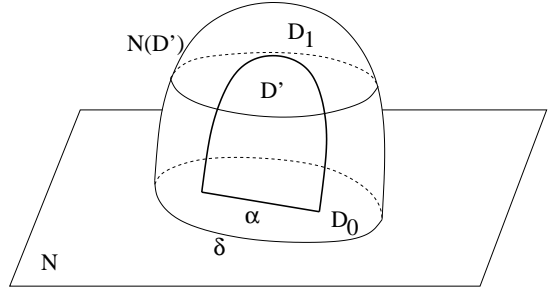}}}
	\caption{Engulfiing $D'$}
	\label{fig17}
\end{figure}

Consider an outermost arc $\alpha$ of $D\cap N$ on $D$.  Then  $\alpha$,
together with an arc $\beta$ of $\bdry D$,
bounds a half-disk $D'\subset
D$ which has no other intersection with $N$.    Although $t(\alpha,Fr_{D'})\leq -1$, it is possible
that $t(\beta, Fr_{D'})>0$, and $D'$ cannot necessarily be made convex.  We instead engulf
$D'$ inside a small neighborhood $N(D')$ of $D'$ so that $N(D')\cap N$ is a disk $D_0$
with boundary $\delta$, $D'\cap D_0=\alpha$, and $\bdry N(D')=D_0\cup_{\delta} D_1$.
See Figure \ref{fig17}.
Use the Legendrian realization principle to make $\delta$ Legendrian.  Then $D_1$ is a convex
disk with Legendrian boundary.
Therefore, pushing a half-disk $D'$ is now reduced to pushing a 3-ball (with convex boundary)
across $N$.

Let $C$ be a closed curve of $D\cap N$ which is innermost on $D$.  ($C$ is homotopically trivial since
$N$ is incompressible.)
$C$ bounds disks $D_1\subset D$ and $D_0\subset N$.
At this point, since $M$ is irreducible, $D_0\cup D_1$ must bound a 3-ball, and we need to push this
3-ball across.   (In the case $M=M_1\#\cdots \# M_n$ in the maximal connect sum decomposition,
we also have $D_0\cup D_1$ bounding a 3-ball.)

\begin{lemma}      \label{breakdown}
Let $\Sigma$ be  a convex surface and $D_0\subset \Sigma$ a convex subdisk with Legendrian
boundary $\gamma$ and $tb(\gamma)<0$.  If $D_1$ is another convex disk with
Legendrian boundary $\gamma$, $\Sigma\cap D_1=\gamma$, and $D_0\cup D_1$ bounds a
tight 3-ball $B^3$, then there exists a sequence of bypass moves which takes
$\Sigma$ to $\Sigma'=(\Sigma\backslash D_0)\cup D_1$.
\end{lemma}

\begin{proof}   In the proof we will alternate between viewing $D_1$ as being parallel to
$D_0$ near $\gamma$ (called $D^h_1$) and perpendicular to $D_0$ (called $D^v_1$).  We may alternate
between the two viewpoints by making modifications to $D_1$ near $\bdry D_1$.
See Figure \ref{alternate}.
\begin{figure}
	{\epsfysize=1.5in\centerline{\epsfbox{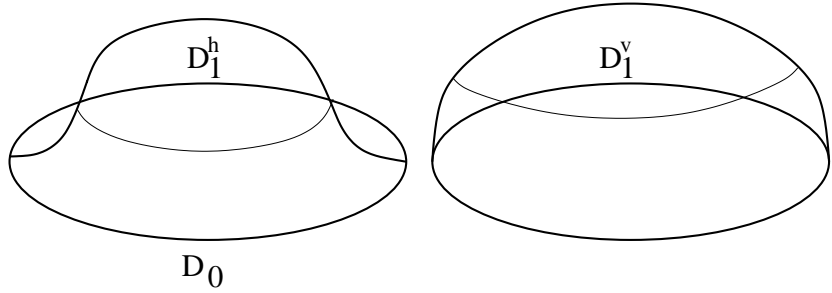}}}
	\caption{Alternating between $D_1^h$ and $D_1^v$}
	\label{alternate}
\end{figure}
Since the endpoints of $\Gamma$ on $D^h_1$ and $D_0$ overlap, this enables us to compare
$\Gamma_{D^h_1}$ and $\Gamma_{D_0}$.  If  $\Gamma_{D^h_1}=\Gamma_{D_0}$, then $\Sigma$
is already contact isotopic to $\Sigma'$.  This follows from the fact that the tight contact structure
on $B^3$ is unique (Theorem \ref{thm:unique}), coupled with the observation that the $I$-invariant contact
structure on $D_0\times I$ is tight.   Otherwise, let $\delta$ be a $\bdry$-parallel dividing curve
on $D_1^v$, and $B_\delta$ be the `corresponding' bypass half-disk which contains $\delta$.
Now, take a parallel copy $(D_1^v)'$ of $D_1^v$ as well as a parallel copy $B'_\delta\subset (D_1^v)'$ of
$B_\delta$.  Attach $B'_\delta$ onto $D_0$.  Then the new $D_0$ agrees with $D_1^h$ along
$B_\delta$ and $\delta$.  If we modify $D_0\mapsto D_0\backslash B'$  and $D_1\mapsto D_1\backslash B$,
then we can induct on $tb(\gamma)$.
\end{proof}

We push the 3-balls across in the following order.  Start with closed curves of
$D\cap N$ which are innermost on $D$, and push the subdisks of $D$ bounding these closed curves
across $N$.     Then continue until $D\cap N$ has no more closed curve components.  Next,
push the half-disk components of $D\backslash N$ which are outermost on $D$.  Continue until
$D\cap N$ is eliminated.  Since each step can be broken down into a sequence of bypass moves by
Lemma \ref{breakdown}, we have extricated $D$ from $N$ through a sequence of bypass
moves attached onto $N$.
This proves the gluing portion of Theorem \ref{classification}.

\subsubsection{Isotopy discretization}
We are now left with the distinguishing process.
We prove that two tight contact structures $\xi$ and $\xi'$ are not isotopic if they do not
correspond to the same connected component of $\mathcal{C}_0$.  For this, we use an ingenious
idea due to Colin \cite{Co97}, called {\it isotopy discretization}.

\begin{lemma}[Isotopy discretization] Let $\xi$ be a tight contact structure on $M$ and $N$, $N'$ be two
convex surfaces with identical Legendrian boundary which are isotopic (but not necessarily
contact isotopic).   Then there exists a sequence of allowable state transitions from
$(\Gamma_N,\xi_{M\backslash N})$ to $(\Gamma_{N'},\xi_{M\backslash N'})$.
\end{lemma}

\begin{proof}  Let $\phi_t:N\rightarrow M$, $t\in[0,1]$, be a 1-parameter family of embeddings,
where $\phi_0$ maps $N$ identically onto $N\subset M$, $\phi_1(N)=N'$, and $\phi_t$ is
independent of $t$  along $\bdry N$.  Break $[0,1]$ into small intervals $I_i=[t_{i-1},t_i]$
with $t_0=0<t_1<\cdots<t_k=1$,
and $k>>0$, so that for each interval $I_i$,
$\phi_{I_i}(N)$ can be sandwiched inside an $N\times I$ with boundary $N^0$ and $N^1$.
Write $N_t=\phi_t(N)$.
The strategy is to compare $N_{t_{i-1}}$ to $N^0$, and to show that we can pass between then
via allowable state transitions, and then compare $N^0$ to $N_{t_i}$.  This approach has the advantage
that every time we are comparing nonintersecting copies of $N$.

Assume therefore that $N$ and $N'$ are nonintersecting (except along their boundary), and
that they bound a layer $N\times I$.
To show that there exists a sequence of bypass moves from $N\times\{0\}$ to $N\times \{1\}$, we
use an important idea due to Giroux \cite{Gi99}, the {\it convex movie}.  According to Giroux,
we may assume (after some perturbations) that there exist $t_i$, $i=1,\cdots, k$, with $0
< t_1 <\cdots < t_k<1$, such that:
\be
\item $N\times \{t\}$ is convex if $t$ does not equal any $t_i$.
\item $N\times \{t_i\}$ is not convex because there is a retrogradient saddle-saddle connection
connecting from a negative hyperbolic singularity to a positive hyperbolic singularity.
\item The saddle-saddle connection serves as a switch as we move from $N\times \{t_i-\varepsilon\}$
to $N\times \{t_i+\varepsilon\}$ ($\varepsilon$ small).
\ee
Moreover, this switching is equivalent to a bypass move.
\end{proof}

This completes the proof of Theorem \ref{classification}.

\subsubsection{Proof of Theorem \ref{gluing}}

Use the same notation as in Section \ref{handlebody}.
Theorem \ref{gluing} is Theorem \ref{classification} rephrased in more combinatorial terms.
First, we may restrict our attention to {\it non-trivial} bypass attachments, which give  state
transitions from one configuration to another.  This follows from the Triviality Lemma (Lemma \ref{triviality}),
since trivial bypass attachments onto a configuration $C$ will only yield an $I$-invariant neighborhood
$D_i\times I$.
Next, we determine when there may exist a state transition $C\sc C'$, or, equivalently, a non-trivial
bypass attachment on some $D_i$.  Draw a Legendrian arc $\delta$ on $D_i$ for which an
attachment of a (theoretical) bypass along $\delta$ from the interior of $B^3=M\backslash (\cup_i D_i)$
still gives $\#\Gamma=1$. The Existence Lemma (Lemma \ref{existence}) then guarantees the actual
existence of a bypass along $\delta$.   Thus, we can peel off a non-trivial $D_i\times I$ layer from
$B^3$.

According to Theorem \ref{classification}, we must reattach the $D_i\times I$ layer `to the other side' (on
$B^3$) to get from $C$ to $C'$.
The state transition $C\sc C'$ is {\it allowable}, if $\Gamma_{\bdry B^3}$ does not change under this
attachment.  The Triviality Lemma again asserts that we are simply attaching a trivial $S^2\times I$
layer, if $\Gamma$ does not change.
Therefore, condition (3) of an
allowable state transition in Theorem \ref{gluing} is equivalent to condition (3) in Theorem
\ref{classification}.  This proves
Theorem \ref{gluing}.

Note that Theorem \ref{gluing} gives a finite, purely combinatorial condition for the contact
structure on the handlebody to be tight.  We are currently working on a computer implementation
of this algorithm with Tanya Cofer, a graduate student at the University of Georgia.

\section{Examples}
Here we present two applications of the Gluing/Classification Theorems.  For a further application,
refer to \cite{H3}, where a version of the Gluing Theorem is presented.

\subsection{Legendrian surgery}
In this section we prove the following result:

\begin{thm} \label{leg} There exists a handlebody $M$ of genus $g=4$ with a tight contact structure $\xi$
which becomes overtwisted after a particular Legendrian surgery.
\end{thm}

This example therefore answers in the negative
the following question:

\begin{q}
Let $(M,\xi)$ be a tight contact manifold, and $L$ a Legendrian curve in $M$.  If a Legendrian
surgery is performed on $M$, then is the resulting contact structure tight?
\end{q}

For {\it holomorphically fillable} contact structures (which are tight by a theorem of
Gromov \cite{Gr} and Eliashberg \cite{E91}), Eliashberg proved \cite{E90} that Legendrian
surgery on a holomorphically fillable contact structure yields a holomorphically
fillable contact structure. Y. Eliashberg and V. Colin informed me that it is easy to extend this
result to show that  Legendrian surgery on a {\it weakly symplectically
semi-fillable} contact structure yields a weakly symplectically semi-fillable structure.
This implies the following corollary:

\begin{cor} There exists a handlebody of genus $4$ with a tight contact structure
which cannot be contact-embedded inside any 3-manifold with a weakly symplectically
semi-fillable
contact structure.
\end{cor}

This contrasts with the genus 1 (solid torus) case, where all the tight contact structures
are contact-embedded inside a lens space with a holomorphically fillable structure
\cite{H1,Gi99}.   It also is a good indication that the classification of tight contact structures
on handlebodies is more subtle than the solid torus case.

We now give a definition of Legendrian surgery for a Legendrian curve $L$ inside $(M,\xi)$
a tight contact 3-manifold.  In general, $L$ may not be homologically trivial, so we must
refer to the twisting number $t(L)$ with respect to some framing, instead of the Thurston-Bennequin
invariant $tb(L)$.  Fix a framing so that $t(L)=0$.  Then a standard neighborhood $N(L)=S^1\times
D^2$ of $L$
has convex boundary $\bdry (N(L))$ with $\#\Gamma_{\bdry (N(L))}=2$
and dividing curves of slope $\infty$.  Here we use the convention that the meridional slope is $0$.
We identify (in a slightly
nonstandard manner) $\psi:\bdry (N(L))\rightarrow \R^2/\Z^2$,
where $\pm(1,0)^T$ corresponds to the meridional
direction and $\pm(0,1)^T$ is the direction of the dividing curves.  Also identify
$\psi':-\bdry (M\backslash N(L))=\bdry(N(L))\rightarrow \R^2/\Z^2$.
We now perform $-1$ surgery with respect to this framing.  Let $M'=(M\backslash N(L))\cup_\phi
N(L)$, where $N(L)$ is glued back via $\phi:\bdry (N(L))\rightarrow - \bdry (M\backslash N(L))$
given by
$\left(\begin{array}{cc} 1 & 0\\ -1 & 1\end{array}\right)$.  Here we are identifying
$\bdry (N(L))=\R^2/\Z^2$ using $\psi$ and $- \bdry (M\backslash N(L))=\R^2/\Z^2$ by
using $\psi'$.
Since the dividing sets on
$-\bdry (M\backslash N(L))$ and $\phi(\bdry(N(L)))$ are identical (although the
characteristic foliations may not exactly line up), we may perturb using Giroux's Flexibility
Theorem to perform the gluing.  The procedure just described is exactly the same
as the procedure which is usually called `Legendrian $tb-1$ surgery' in the context of
holomorphically fillable structures.

\begin{figure}
	{\epsfysize=2in\centerline{\epsfbox{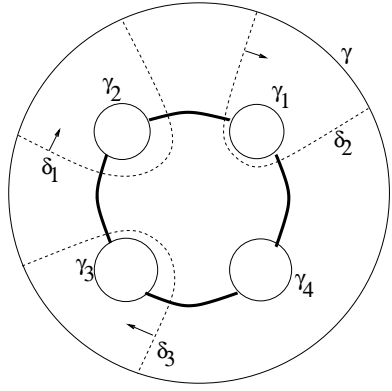}}}
	\caption{Dividing curves on $\Sigma$.  Dark solid lines are dividing curves.}
	\label{fig10}
\end{figure}

\begin{figure}
	{\epsfysize=2in\centerline{\epsfbox{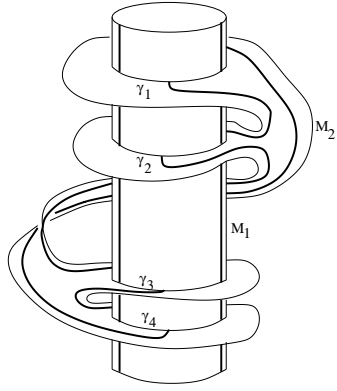}}}
	\caption{Construction of $M$.  The meridian of $M_1$ has slope $+1$.}
	\label{fig11}
\end{figure}

Let us now describe our example.  Take a solid torus $M_1=S^1\times D^2$ and identify the
boundary $T=\R^2/\Z^2$ so that the meridional slope is $+1$.
We take $\Gamma_T$ to consist of two dividing curves of slope $\infty$.  Also
assume the actual characteristic foliation is a ruling by Legendrian curves of
slope $0$ -- these are called {\it Legendrian rulings}.  Then there exists
a unique tight contact structure on $M_1$ with this boundary condition, and
the unique tight structure can be realized as a standard neighborhood of a Legendrian curve $L$.
(See the `Basic building blocks' section of \cite{H1} for more details.)
Next, let $\Sigma$ be a 4-holed convex disk with Legendrian boundary
$\bdry \Sigma=\gamma-(\gamma_1+\gamma_2+
\gamma_3+\gamma_4)$.  Let $t(\gamma)=0$, $t(\gamma_i)=-1$, and $\Gamma_\Sigma$
consist of 4 arcs -- from $\gamma_i$ to $\gamma_{i+1}$, where $i=1,2,3,4$ and
$\gamma_i=\gamma_{i+4}$.
Extend $\Sigma$ slightly beyond the boundary Legendrian curve (still call this $\Sigma$).
Let $M_2$ be an $I$-invariant neighborhood of $\Sigma$.
We let $M=M_1\cup M_2$, where the $\gamma_i$ (as well as parallel copies of $\gamma_i$
on $\bdry M_2$) are attached  to Legendrian ruling
curves of $M_1$.  For two of the attachments, say along $\gamma_1$ and $\gamma_2$,
the normal direction to $\Sigma$ is the
same as that given by $(0,1)^T$, and for the other two, $\gamma_3$ and $\gamma_4$, the
normal direction to $\Sigma$ is the opposite to the one given by $(0,1)^T$.  Refer to
Figures \ref{fig10} and \ref{fig11}.

\begin{lemma}  Let $L$ be a Legendrian core curve of $M_1$  with highest
twisting number inside $M_1$. If $M'$ is obtained from $M$ by Legendrian $-1$
surgery along $L$, then the contact structure on $M'$ is overtwisted.
\end{lemma}

The highest twisting number condition is the same as requiring $M_1\backslash N(L)$
to be an $I$-invariant $T^2\times I$.

\begin{proof}  This is because the slopes of the new meridional disks in $M'$ are $0$, and we
can patch $4$ copies of the meridional disk onto $\Sigma$ to make it into an overtwisted disk.
We are essentially unlinking $\Sigma$ from $M_1$ after surgery.
\end{proof}

\begin{thm}  The contact structure on $M$, described above, is tight.
\end{thm}

Intuitively speaking, we cannot unlink $\Sigma$, which wants to be an overtwisted
disk.

\begin{proof} We apply a variant of the Gluing Theorem.  Instead of cutting along $4$ disks,
we will cut along $3$ disks $D_i=\delta_i\times I$, $i=1,2,3$, where $\delta_i$ are arcs on $\Sigma$
drawn in dotted lines in Figure \ref{fig10}.  Orient $D_i$, where the normal orientation is indicated by the
arrows in Figure \ref{fig10}.
After cutting, $M\backslash (D_1\cup D_2\cup D_3)$
is a solid torus, and its boundary will have two copies each of $D_i$, which we denote $D_i^+$ and $D_i^-$,
depending on whether the orientation induced from $D_i$ agrees or disagrees with the
boundary orientation of $M\backslash (D_1\cup D_2\cup D_3)$.
In order to prove tightness, we will need to use our knowledge about tight contact structures on
solid tori, in exchange for simplifying
the combinatorics somewhat.   (Refer to \cite{H1} for a discussion of tight contact structures on
solid tori.)

\vskip.12in
\noindent
{\bf Initial configuration $C$:}   For each of the $D_i^+$, the initial configuration of dividing curves
is as in Figure
\ref{fig13}. (The portions with $D_i^-$ are not pictured here.)
This is due to the $I$-invariance of the tight contact structure on $\Sigma$.
\begin{figure}
	{\epsfysize=2in\centerline{\epsfbox{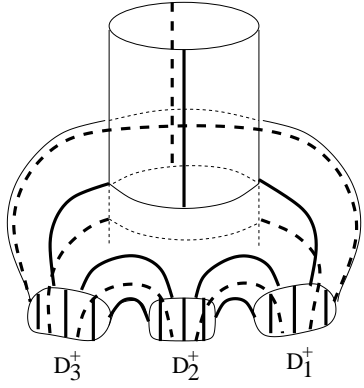}}}
	\caption{Initial configuration}
	\label{fig13}
\end{figure}
If we cut $M$ along $D_1$, $D_2$, $D_3$ and rounded the edges, the resulting solid torus
will  have slope $\infty$, which is isomorphic to $M_1$.
A bypass which is attached to $D_i$ `from the outside' is attached from the direction of the oriented normal
of $D_i$,
and a bypass which is attached `from the inside' is attached from the opposite direction.
For each $D_i^+$, there are three possible dividing curve configurations, given in Figure \ref{fig14} as
$e=a^3$, $a$, and $a^2$ -- they are denoted in suggestive group-theoretic notation to indicate
that each outer bypass attachment is an action by $a$, and three bypasses in a row gives back the
original configuration.  We then write the initial configuration $C$ as $(e,e,e)$ (the $i$-th component
is for $D_i^+$).
\begin{figure}
	{\epsfysize=1.3in\centerline{\epsfbox{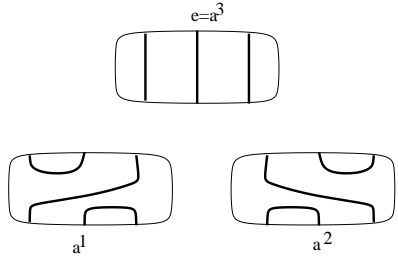}}}
	\caption{Possible configurations for $D_i^+$}
	\label{fig14}
\end{figure}

\vskip.12in
\noindent
{\bf Configuration changes from $C=(e,e,e)$:}
(1) There can be no (nontrivial) bypasses attached to any $D_i^+$
from the inside --
any such bypasses would create a homotopically trivial curve.
This is an easy check once we attach a bypass from the inside and
round the edges in Figure \ref{fig13}.   Therefore, there is no
$e\sc a^2$.

(2) There can exist bypasses from the outside.    Attaching a bypass to
the outside of Figure \ref{fig13} along $D_i^+$  does not change the dividing set of the portion
pictured in Figure \ref{fig13}, after rounding.    $D_i^-$  on $M\backslash (D_1\cup D_2\cup
D_3)$
is given in Figure \ref{fig15}.  Peeling off the bypass on the `outside' (= bypass `on the
inside' of Figure \ref{fig15}) adds an extra twist to the dividing curves of $\bdry(M\backslash
(D_1\cup D_2\cup D_3))$, so
the slopes change from $\infty$ to $-1$.

This proves that $(e,e,e)\sc (a,e,e)$, $(e,e,e)\sc (e,a,e)$, $(e,e,e)\sc (e,e,a)$ are the only allowed
state transitions from $(e,e,e)$.

\vskip.12in
\noindent
Next assume we have a configuration $C'$ where the entries $c_i$ are either $e$ or $a$,
and the set of indices $i$ with $c_i=a$ is a subset of either $\{1,2\}$ or $\{3\}$. In other
words, assume $C'\in {\cal C}' =\{(e,e,e), (a,e,e), (e,a,e), (e,e,a), (a,a,e)\}$.
\vskip.12in
\noindent
{\bf Configuration changes from ${\cal C}'$:}  We claim that if
$C'\in {\cal C}'$, then any allowable state transition will remain in ${\cal C}'$.
We are trying to attach bypasses onto $D_i$ from the inside or from the outside.

(1) If $c_i=a$,
then it is possible to find a bypass on the inside, and use
it to perform a state change $a\sc e$.    This is precisely the inverse process of $e\sc a$, where the
bypass is attached from the outside.

(2)  If $c_i=a$, it is not possible to attach a second bypass from the outside.
Peeling off a second bypass from $M_1$ corresponds to reducing the boundary slope
from $-{1\over n}$ ($n$ is the number of $a$'s in $C'$) to $0$.
However, the sign of the bypass is opposite that of the first bypass, and
it is not possible to find two bypasses of opposite signs on a basic $T^2\times I$
layer from slope $\infty$ to slope $0$.  Recall from \cite{H1} that a basic slice $(T^2\times I,\xi)$
with boundary slopes $s_1=s(\Gamma_{T^2\times \{1\}})=\infty$  and
$s_0=s(\Gamma_{T^2\times \{0\}})=0$ has {\it relative Euler class} $e(\xi,s)=\pm(1,-1)$ (for a
suitable section $s$ on the boundary).  If
$T^2\times I$ admits a splitting $(T^2\times [0,{1\over 2}])\cup (T^2\times [{1\over 2},1])$ with
$s_{1\over 2}=-1$,
and $T^2\times [0,{1\over 2}]$  has relative Euler class $e(\xi,s)=(0,1)$, then
the relative Euler class for $T^2\times [{1\over 2},0]$ is dictated to be $(-1,0)$.  A second bypass
from the outside violates this consistency condition.
(Another way to obtain a contradiction is to observe that
stacking a $D^2\times [1,2]$ layer atop $D^2\times [0,1]$, with $\Gamma_{D^2\times \{i\}}=a^i$,
$i=0,1,2$, gives an overtwisted disk on $\bdry (D^2\times [0,2])$.)

(3)  Assume $c_i=e$.  If $C'\not= (a,a,e)$, then at least two $c_j=e$, and
a bypass on the inside gives rise to an overtwisted disk. If $C'=(a,a,e)$, then the configuration in the
neighborhood of $D_1^+\cup D_2^+\cup D_3^+$,
after smoothing $\bdry D_1^+$ and $\bdry D_2^+$,
is equivalent to the configuration in Figure  \ref{fig15}.
\begin{figure}
	{\epsfysize=1.7in\centerline{\epsfbox{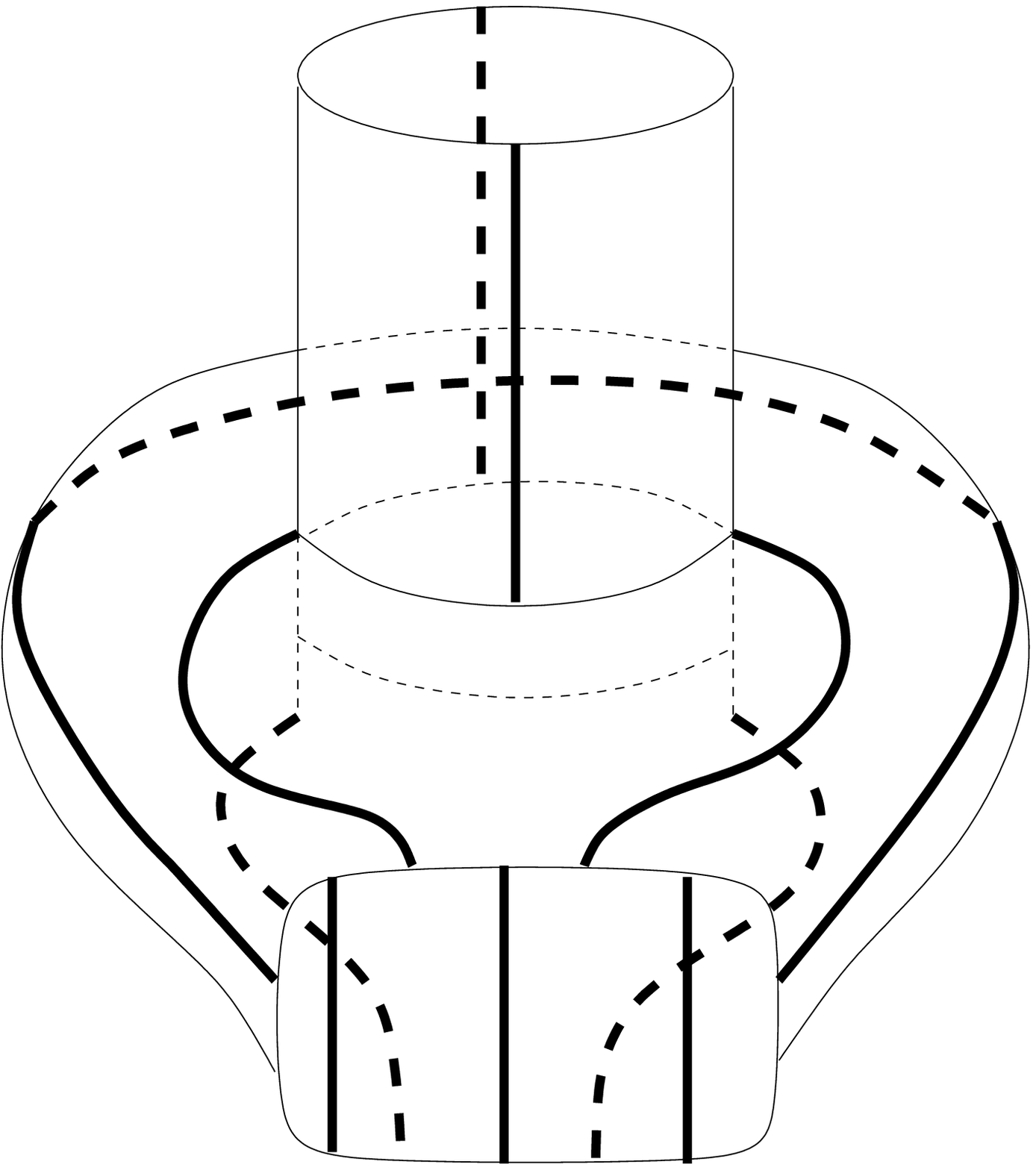}}}
	\caption{}
	\label{fig15}
\end{figure}
There cannot exist a bypass from the inside for this case, because there will be two bypasses
of opposite signs on a basic $T^2\times I$ as in (2).

(4)  If $c_i=e$, then it is possible to find a bypass on the outside, but only if $i=1$ or $2$, and $C'=
(a,e,e)$ or $(e,a,e)$.
This follows from the same reason as (2).

\vskip.12in
\noindent
All of the cut-open solid tori are
then tight.\end{proof}

This construction can be generalized to any $g\geq 4$, as long we have at least two
`positive' attachments of $\gamma_i$ onto $M_1$ and at least two `negative' attachments.
It also turns out that this configuration -- a Legendrian curve $L$ with a candidate overtwisted
disk $\Sigma$ linked onto $L$ -- is the only obstruction to proving that Legendrian
$tb-1$ surgery on a given tight contact manifold produces a tight contact manifold.
The following question is therefore
equivalent to the existence of such configurations:

\begin{q}
Does every Legendrian $-1$ surgery on a closed tight contact manifold produce a
tight contact manifold?
\end{q}

\subsection{Dehn surgeries on $T^3$}

Start with $(T^3,\xi_n)$, where $T^3=S^1\times T^2=\R^3/\Z^3$ with coordinates $(z,x,y)$, and
$\xi_n$ is given by the standard 1-form $\alpha_n=\sin (2\pi n z)dx+\cos (2\pi n z)dy$ with $n$
complete twists.
Take $\phi_{k_1,k_2}:T^3\rightarrow T^3$ given by $(z,x,y)\mapsto (k_1x+k_2y+z,x,y)$, where
$k_1, k_2\in \Z$.  Then set $\xi_{n,k_1,k_2}=\phi_{k_1,k_2}^*\xi_n$.  These are all nonisotopic
on $T^3$, due to Kanda \cite{K97} and Giroux \cite{Gi99}.  (It would also follow
from the proof of Theorem \ref{T^3} below.) The integers $k_1$, $k_2$ are usually called
{\it holonomy}.
Let $L$
be a Legendrian curve isotopic to the fiber $S^1$ with $t(L)=-n$.    This
twisting number is measured relative to any $T^2$, containing $L$, which is
isotopic to $S^1\times \{x=0\}$.  (Note that this twisting number is independent of the
base curve, although $\{x=0\}$ was used.)  Consider a standard neighborhood $N(L)$,
and identify (in a nonstandard manner)
$-\bdry (T^3\backslash N(L))=\bdry (N(L))=\R^2/\Z^2$ by letting
$\pm (1,0)^T$ be the meridional direction and $\pm (0,1)^T$ be the preferred longitude arising
from the fiber direction.

We study $r$-Dehn surgeries, where
$-\infty<r <-n$.  The meridional slopes $s$ under the identification with $\R^2/\Z^2$
satisfy $-{1\over n}< s={1\over r}<0$.  Notice that $\bdry (N(L))$ has slope $-{1\over n}$, and
we simply fill in with any tight solid torus $(S^1\times D^2,\zeta)$ with boundary slope $s$,
$\#\Gamma_{\bdry (S^1\times D^2)}=2$, and
meridional slope ${1\over r}<0$.     Observe that $\zeta$ may be {\it virtually overtwisted}.
(Recall a tight contact structure is {\it virtually overtwisted} if it becomes overtwisted after a lift to
some finite cover.  This contrasts with the notion of a {\it universally tight} contact structure, which remains
tight even after passage to the universal cover.)
Write $M$ for this Seifert fibered space over $T^2$ with
invariant $-s$,
and $\eta(n,k_1,k_2,\zeta)$ for the glued contact structure.   Notice that if we fix $r$ and
$k_1, k_2$, then there are finitely many
$\eta(n,k_1,k_2,\zeta)$'s.   However, if we only fix $r$, then there exist infinitely many
$\eta(n,k_1,k_2,\zeta)$'s.

\begin{thm}  \label{T^3} The contact structures $\eta(n,k_1,k_2,\zeta)$ are tight, and
$\eta(n,k_1,k_2,\zeta)$ and
$\eta(n',k_1',k_2',\zeta')$ are isotopic if and only if $n=n'$, $k_i=k_i'$, ($i=1,2$), and $\zeta=\zeta'$.
\end{thm}

If we took $\zeta$ to be virtually overtwisted, then we have the following corollary.

\begin{cor} There exist closed toroidal 3-manifolds which have infinitely many virtually overtwisted
tight contact structures (up to isotopy).
\end{cor}

Theorem \ref{T^3} generalizes easily to Seifert fibered spaces over $T^2$ with multiple singular
fibers and Seifert invariants $(-s_1,-s_2,\cdots,-s_k)$.

\begin{proof}[Proof of Theorem \ref{T^3}]
Consider the following convex decomposition  $M\stackrel{S_1}\sa M'\stackrel{S_2}\sa M''=S^1\times
D^2$, where $S_1=S^1\times \{y=0\}$, and $S_2=S^1\times \{x=0\}$. Properly speaking,
$S_2$ is an annulus since the
first cutting has already taken place in the decomposition.
We prove the theorem by tracing back along the convex decomposition.    Note that we will often make
modifications using the Flexibility Theorem or the Legendrian realization principle without explicitly stating
that they will be used.

\s
\noindent
{\bf Initial Configuration $C_0$.}  Initially, $\Gamma_{S_1}$ will consist of $2n$ parallel
dividing curves with holonomy $k_1$ and $\Gamma_{S_2}$ will consist of $2n$ parallel
curves with holonomy $k_2$.

\s
\noindent
{\bf Inductive Assumption.}  Suppose we have inductively (via state transitions) at a configuration
$C'$ where:
\be
\item $\Gamma_{S_1'}$ consists of $2(n+m)$, $m\geq 0$, parallel dividing curves with holonomy $k_1$.
\item $S_1'$ is obtained by `folding' inside an $I$-invariant neighborhood $S_1\times I$ of $S_1$. The
folding operation is described in detail in Section 5.3 of \cite{H1}, and is used to increase the number of
parallel dividing curves.
\item $\Gamma_{S_2'}$ consists of $2n$ parallel arcs which go across from one side of the annulus to
the other, together with $2m$ one-sided ($=$ both endpoints on one side of the annulus) arcs.
Half of the one-sided arcs have endpoints on one side of the annulus, and half have endpoints on the other.
Moreover, $\Gamma_{S_2'}$, when glued up along $S_1'$, will become $2n$ parallel curves with
holonomy $k_2$.
\item The tight contact structure on $M\backslash (S_1'\cup S_2')$, after rounding, is $\zeta$.
\ee
Let $\mathcal{C}'$ be the set of such configurations.

\s
\n
{\bf State changes along $\Gamma_{S_2'}$.}  We claim there can be no state transitions.  To do this,
we examine all possible locations on $\Gamma_{S_2'}$ where a non-trivial bypass may be attached.
Let $\delta$ be an arc of attachment for the candidate bypass, $p_i$, $i=1,2,3$, the points of intersection
with $\Gamma_{S_2'}$, and $\gamma_i$ be the dividing curve containing $p_i$.  We enumerate all
the possible cases and eliminate them in turn, by examining $\bdry (M\backslash (S_1'\cup S_2'))$.  The following
key observation helps reduce the potentially infinite number of possibilities to a finite number:

\s
\n
{\it Key observation.}  Any $\bdry$-parallel arc on $S_2'$ which does not intersect $\delta$ may
be pushed into the $S'_1$ portion.  See Figure \ref{push}.
\begin{figure}
	{\epsfysize=1.5in\centerline{\epsfbox{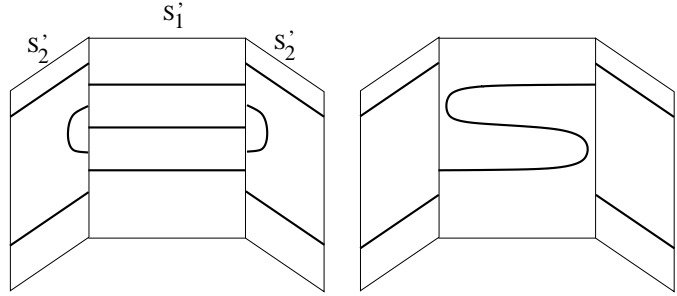}}}
	\caption{This represents only a portion $S_1'$ and $S_2'$. }
	\label{push}
\end{figure}

\s\n
The various cases are enumerated in Figure \ref{cases}. The dark solid lines are the dividing curves
and the light solid lines are the candidate $\delta$. The candidate bypasses may be attached to
the front or the back.  Note that all of the extraneous
dividing curves have already been pushed across.  All the cases fail because the candidate bypass
(attached from the interior of $S^1\times D^2$) gives rise to either (1) a homotopically trivial dividing
curve or (2) a convex torus on the interior of $S^1\times D^2$ with zero boundary slope, which is a contradiction,
since the boundary slope of $S^1\times D^2$ is $-{1\over n}<0$ and the meridional slope is
${1\over r}<0$.   We will treat a few cases, and leave the rest to the reader.
\begin{figure}
	{\epsfysize=5in\centerline{\epsfbox{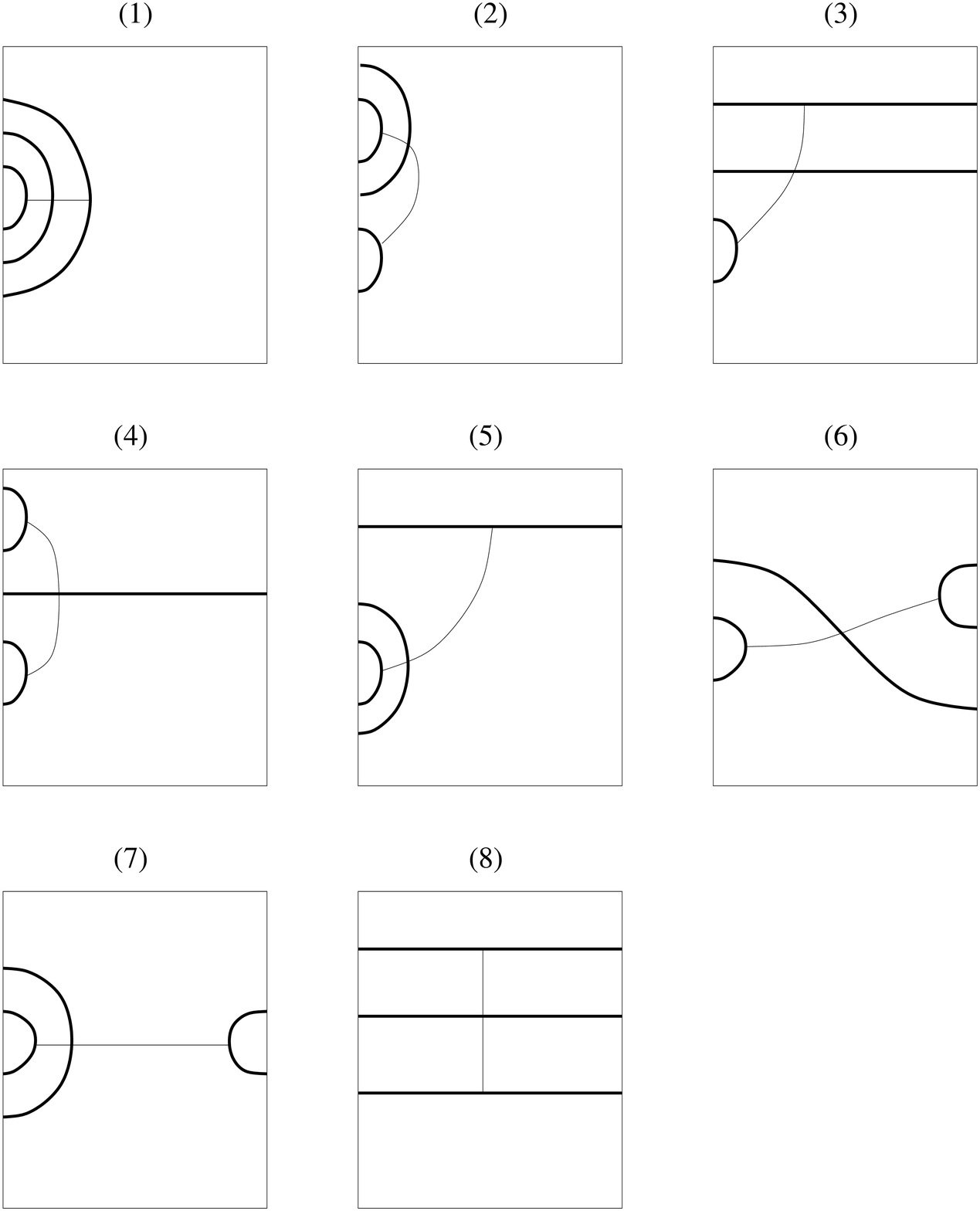}}}
	\caption{}
	\label{cases}
\end{figure}

\s\n
{\bf Case 1.}  See Figure \ref{case1}.    Here, regardless of which side we attach a bypass along
$\delta$, there will exist a homotopically trivial disk after attachment.  In the figure only the relevent portions
of $S_i'$ are shown, and the edges have already been rounded.  The two thin lines are the arc of attachment
$\delta$, and the hypothetical bypasses are attached from the front onto either of the two thin lines.
\begin{figure}
	{\epsfysize=1.5in\centerline{\epsfbox{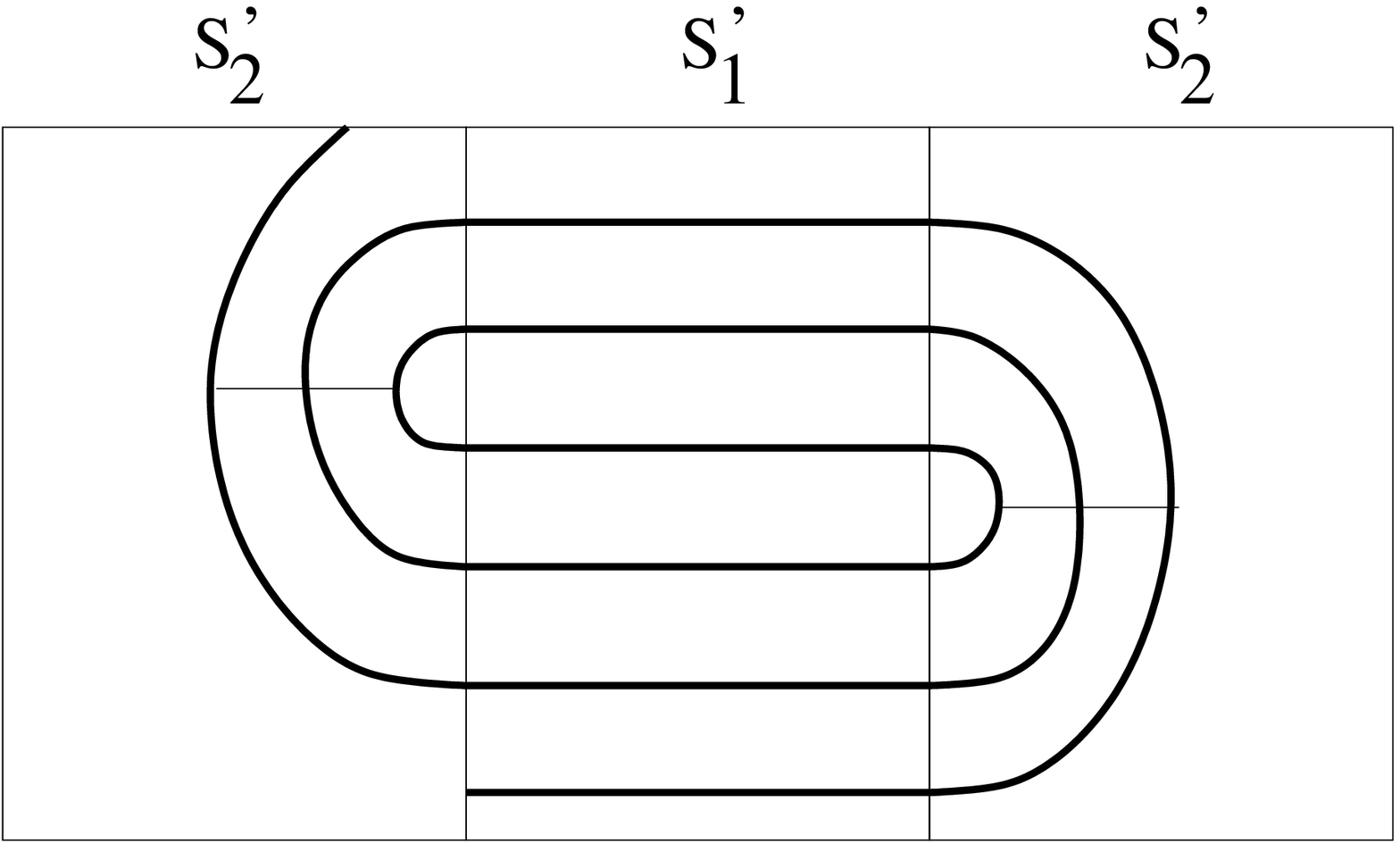}}}
	\caption{Case 1}
	\label{case1}
\end{figure}

\s\n
{\bf Case 3.}  See Figure \ref{case3}.
\begin{figure}
	{\epsfysize=1.5in\centerline{\epsfbox{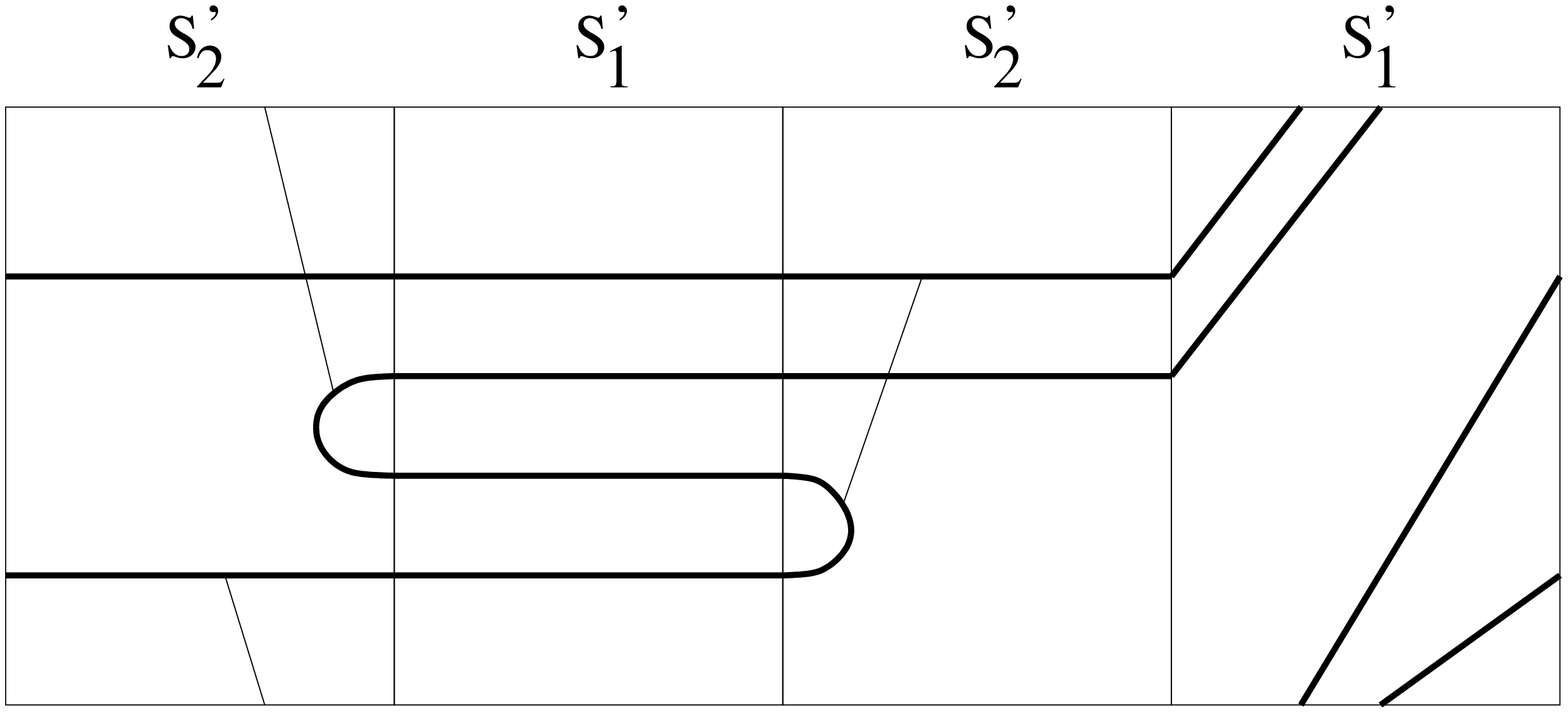}}}
	\caption{Case 3}
	\label{case3}
\end{figure}
Here we have depicted the case $n=1$.  The left and right sides of the large rectangle are  identified
to give $\bdry (S^1\times D^2)$. The $\delta$ to the right immediately gives a homotopically trivial
dividing curve after attachment.  The $\delta$ to the left yields a dividing set with slope zero.
This gives us a contradiction.

\s\n
{\bf Case 6.} See Figure \ref{case6}.
\begin{figure}
	{\epsfysize=1.5in\centerline{\epsfbox{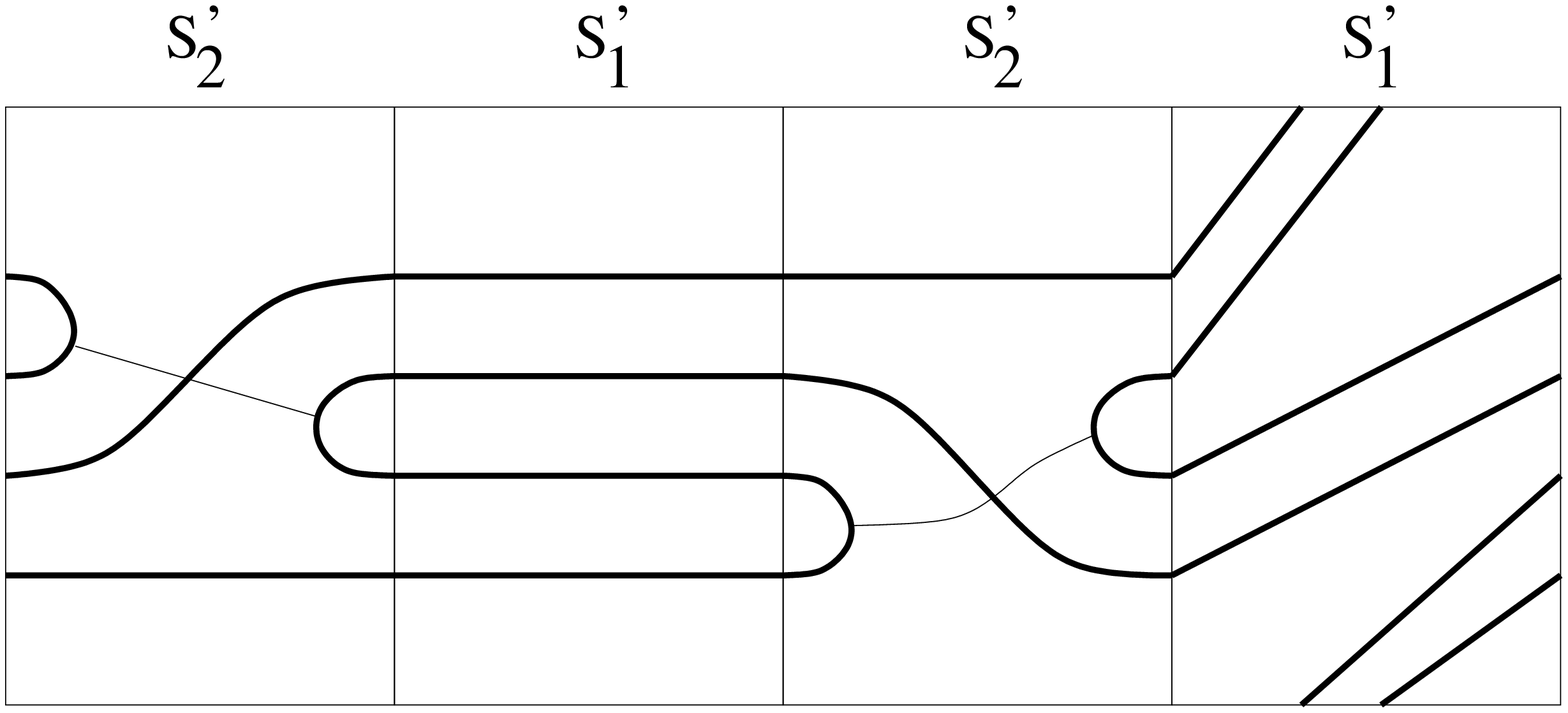}}}
	\caption{Case 6}
	\label{case6}
\end{figure}
This is similar to Case 3.  The $\delta$ to the right gives a homotopically trivial curve, and
the $\delta$ to the left yields zero slope.

\s\n
{\bf Case 8.}  In this case we may use the key observation to reduce to the case when each of the
$S_i'$ has exactly $2n$ dividing curves (all the one-sided components can be pushed across).
The attachment of a bypass along $\delta$ is tantamount to modifying the boundary slope of
$\bdry (S^1\times D^2)$ from $-{1\over n}$ to $0$.

\s\n
Therefore we have proved that there are no state transitions which modify $\Gamma_{S_2'}$.

\s\n
{\bf State changes along $\Gamma_{S_1'}$.}  We show that any state transition still leaves us in
$\mathcal{C}'$.     Assume $S_1'\sc S_1''$.  Here $S_1''$ is  parallel to and disjoint from $S_1'$ in most
cases except when $n=1$.  (See Case 3 below for the exception.)
Let $N_1$ be the $I$-bundle bounded by $S_1'$ and $S_1''$.
We have three cases.

\s\n
{\bf Case 1.}  $\#\Gamma_{S_1''}=\#\Gamma_{S_1'}+2=2(n+m)+2$.  Let $\gamma'$ be a vertical Legendrian
curve on $S_1''$ with $t(\gamma')=-(n+m)$, and let $\gamma''$ be a vertical Legendrian curve on
$S_1''$ with $t(\gamma'')=-(n+m)-1$.  Let $A$ be an annulus with $\bdry A=\gamma'\cup\gamma''$ and
$A\cap N_1=\bdry A$.  Then, using the Imbalance Principle in \cite{H1}, we find
a bypass along $\gamma''$.  Attaching the bypass is equivalent to finding an $I$-invariant
layer $N_2\subset M\backslash int(N_1)$   with $\bdry N_2=S_1''\cup S_1'''$, where $\Gamma_{S_1'''}=\Gamma_
{S_1'}$.  We find that $N_1\cup N_2$ is therefore an $I$-invariant tight (one-sided) neighborhood of $S_1'$.
Now, since $S_1'$ itself is contained in an $I$-invariant neighborhood of $S_1$,
we conclude the same for $S_1''$.  This proves conditions (1) and (2) of the inductive assumption.
(3) and (4) are immediate from the invariance of $S_2'$ with fixed $S_1'$.

\s\n
{\bf Case 2.}  $\#\Gamma_{S_1''}=\#\Gamma_{S_1'}-2=2(n+m)-2$.  Assume $m>0$ or
$n>1$. Take a vertical Legendrian $\gamma'$ on
$S_1'$ with $t(\gamma')=-(n+m)$ and a vertical Legendrian $\gamma''$ on $S_1''$ with
$t(\gamma'')=-(n+m)+1$.  Take an annulus $B$ with $\bdry B=\gamma'\cup \gamma''$
and $B\subset N_1$.
$B$ can be extended to an annulus $S_2''$ isotopic to $S_2'$.  However, since no state transitions can
occur on $S_2'$, we may assume that $S_2''=S_2'$.  Moreover, we may assume that $S_1''$ is obtained
from $S_1'$ by attaching a bypass which lies on $S_2'$.  Now, take
all the one-sided dividing arcs on $S_2'$ which end on a fixed boundary component of $S_2'$.
Attach all the bypasses corresponding to these one-sided dividing arcs to get an isotopic copy of
$S_1$.  Do the same with the other boundary component.  Then we get another isotopic copy of
$S_1$, and $S_1''$ is sandwiched inbetween.
This proves the conditions of the inductive assumption.

\s\n
{\bf Case 3.}  We have one more case left, if $\#\Gamma_{S_1'}=2$ ($n=1$).  Without loss of generality,
assume that $k_1=0$.  Then the slope of $\Gamma_{S_1'}$ is $0$.  A non-trivial bypass attachment will
give $S_1''$ with slope $-{1\over k}$, $k\in \Z$.  If we consider the tight contact structure on $N_1$, then,
by the classification of tight contact structures on $T^2\times I$ (c.f. \cite{H1}),
we may take $k\in \Z$ to be large positive.
Now consider vertical Legendrians $\gamma'$ on $S_1'$ with $t(\gamma')=-1$ and
$\gamma''$ on $S_1''$ with $t(\gamma')=-k$.  Take an annulus $A$ with $\bdry A=\gamma'\cup\gamma''$
and $A\cap N_1=\bdry A$.
Then there will exist bypasses along $\gamma''$ which
allow us to find $S_1'''$ with slope $-1$.  Now, cutting along $S_2'$ (and using the
fact that $\Gamma_{S_2'}$ does not depend on the cut), we see that the bypass attachment
corresponding to $S_1'\sc S_1''$ will yield a slope $0$ solid torus inside the
slope $-1$ solid torus $M\backslash (S_1'\cup S_2')=S^1\times D^2$, a contradiction.
A similar but easier argument shows that there cannot be a dividing curve decrease if $n>1$ and $m=0$.
\end{proof}

\vskip.2in
\noindent
{\it Acknowledgements and notes:}  I would like to thank John Etnyre for many ideas and suggestions.
I would also like to thank Will Kazez and Gordana Mati\'c for their
enthusiasm and suggestions during our informal seminar at the University of Georgia, where
I first presented most of this material.  In particular, the formulation of the gluing theorem, although
still cumbersome, has been much simplified from the original version, thanks to Gordana Mati\'c.
Vincent Colin recently obtained an example of a tight contact structure which becomes overtwisted
after an admissible transverse surgery (see \cite{Co00}).  The example of Theorem \ref{leg} (minus the
proof of tightness) was also known to him.


\end{document}